\documentclass[a4paper]{article}

\usepackage{amsthm}
\usepackage{amssymb}
\usepackage{amsmath}
\usepackage[all]{xy}
\usepackage{graphicx}

\begin{document}

\baselineskip15pt
\textwidth=12truecm
\textheight=20truecm
\hoffset-.1cm
\voffset-.5cm
\font\eightrm=cmr8


\newtheorem{thm}{Theorem}[section]
\newtheorem{cor}[thm]{Corollary}
\newtheorem{lem}[thm]{Lemma}
\newtheorem{prop}[thm]{Proposition}
\newtheorem{defn}[thm]{Definition}
\newtheorem{ax}{Axiom}

\theoremstyle{definition}
\newtheorem{ex}{Example}[section]

\theoremstyle{remark}
\newtheorem{rem}{Remark}[section]
\newtheorem*{notation}{Notation}

\numberwithin{equation}{section}

\newcommand{\thmref}[1]{Theorem~\ref{#1}}
\newcommand{\secref}[1]{\S\ref{#1}}
\newcommand{\lemref}[1]{Lemma~\ref{#1}}
\newcommand{\propref}[1]{Proposition~\ref{#1}}
\newcommand{\defnref}[1]{Definition~\ref{#1}}
\newcommand{\bysame}{\mbox{\rule{3em}{.4pt}}\,}

\newcommand{\mA}{A}
\newcommand{\mB}{B}
\newcommand{\mC}{\mathcal C}
\newcommand{\mE}{\mathcal E}
\newcommand{\mF}{\mathcal F}
\newcommand{\mL}{\mathcal L}
\newcommand{\mM}{\mathcal M}
\newcommand{\mO}{\mathcal O}
\newcommand{\mX}{\mathcal X}
\newcommand{\mY}{\mathcal Y}
\newcommand{\mZ}{\mathcal Z}

\newcommand{\cM}{{\mathcal M}{\it {od}}}
\newcommand{\cB}{{\mathcal B}{\it {imod}}}

\newcommand{\sC}{{\it C*}-}

\newcommand{\bC}{{\mathbb C}}
\newcommand{\bN}{{\mathbb N}}
\newcommand{\bT}{{\mathbb T}}
\newcommand{\bZ}{{\mathbb Z}}
\newcommand{\bM}{{\mathbb M}}

\newcommand{\coe}{\mO_{\mE}}
\newcommand{\oex}{\mO_e}
\newcommand{\oro}{\mO_\rho}

\newcommand{\om}{\omega}

\newcommand{\ra}{\rightarrow}
\newcommand{\hra}{\hookrightarrow}

\author{Ezio Vasselli \thanks{Dipartimento di Matematica dell'Universit\`a
di Roma ``Tor~Vergata''. E-mail {\tt vasselli@mat.uniroma2.it}}}

\title{Continuous Fields of \sC Algebras Arising from Extensions of Tensor \sC Categories}

\maketitle

\abstract{The notion of extension of a given \sC category $\mC$ by a \sC algebra $\mA$ is introduced. In the commutative case $\mA = C(\Omega)$, the objects of the extension category are interpreted as fiber bundles over $\Omega$ of objects belonging to the initial category. It is shown that the Doplicher-Roberts algebra (DR-algebra in the following) associated to an object in the extension of a strict tensor \sC category is a continuous field of DR-algebras coming from the initial one. In the case of the category of the hermitian vector bundles over $\Omega$ the general result implies that the DR-algebra of a vector bundle is a continuous field of Cuntz algebras. Some applications to Pimsner \sC algebras are given.}

\section{Introduction.}

It is well known that tensor categories play an important role in the theory of duality. In particular in the basic paper \cite{DR89} suitable abstract strict tensor \sC categories are characterized as duals of compact groups. One of the steps of the proof of that result is the association of a \sC algebra, called the DR-algebra, to each object of the given strict tensor \sC category. Some of these \sC algebras have been studied (\cite {DR87, CP92, PR96, DPZ97, KPW98}) and others are well known (\cite{Cun77}); the aim of this work is to start the study of DR-algebras arising from a certain class of strict tensor \sC categories, having not trivial space of intertwiners of the identity object. Such categories have already been studied in the context of superselection structures, in the case where the fixed point algebra of a crossed product has a non trivial center (see \cite {BL97} and related references). The present paper is organized as follows:

In section 2 some basics are given about strict tensor \sC categories, DR-algebras and continuous fields of \sC algebras with a circle action. Furthermore some preliminary results are illustrated.

In section 3 we introduce a simple procedure for extending a given \sC category $\mC$ by a unital \sC algebra $\mA$. The objects of the extended category $\mC^\mA$ are projections in the spatial tensor product $\mA \otimes (\rho,\rho)$ for some object $\rho$ in $\mC$. We give in the commutative case $\mA = C(\Omega)$ a picture of the extension as the categorical analogue of a fiber bundle over $\Omega$, and show that it inherits naturally the tensor structure by $\mC$. We also give a description of $\mC^{C(\Omega)}$ in terms of principal bundles, associating to objects of the extension category elements of suitable cohomology sets (Prop. \ref{thm_co}). 

In section 4 we study DR-algebras associated to objects of $\mC^\mA$. Our main result is that in the commutative case they have the structure of continuous fields of \sC algebras, whose fiber is given by DR-algebras associated to objects of $\mC$ (Thm. \ref{main-thm}). Particular care is given at the case of Hermitian vector bundles, whose associated DR-algebra is a continuous field of Cuntz algebras. We prove that the isomorphism between the associated DR-algebras does not imply the isomorphism between the vector bundles (Prop. \ref{prop_iso}), so that classification questions arise.

In section 5 we return to the case in which $\mA$ is non commutative and investigate the existence of tensor products on our extension categories. Further structure has to be added to obtain a tensor product: we need an abstract analogue of the left action on a Hilbert \sC bimodule, so we consider endomorphisms of the form $\mA \ra \mA \otimes (\rho,\rho)$ (see \cite {DR89}, \S 1.7). With a technical assumption (namely the existence of a unique \sC norm over $(\rho,\rho) \otimes_{alg} \mA$ extending the \sC norms of the factors for each object $\rho$ in $\mC$) we construct a tensor product on the extension category (Prop.\ref{prop_tens}). As an application of the abstract results we characterize Hilbert \sC bimodules of the form $\mX = \mE \otimes_{\mZ} \mA$, where $\mZ$ denotes the center of $\mA$ and $\mE$ is the module of sections of a vector bundle over $spect(\mZ)$, in terms of properties of the Pimsner \sC algebra associated to $\mX$ (Prop. \ref{prop_mod}).

Some of the results in section 4 have been proved independently (and previously) by J.E. Roberts \cite {Rob}.

\section{Preliminaries.}

\subsection{DR-Algebras.}

A detailed exposition of tensor categories can be found in \cite{McL}; in \cite{DR89} strict tensor \sC categories with conjugates are considered and illustrated by examples, and the notion of DR-algebra is introduced (in both references the term {\em monoidal} is used instead of {\em tensor}). More recent developements on strict tensor \sC categories can be found in \cite{LR97}; we refer to this paper for the definition of conjugates.

To each object $\rho $ of a strict tensor \sC category $\mC$ is canonically associated a \sC algebra $\oro $, whose construction we briefly recall. We consider for each $k$ in ${\bZ}$ the Banach space $\oro^k$ defined by the inductive limit

$$
...\stackrel {\times 1_\rho} \longrightarrow
( \rho^r,\rho^{r+k}) \stackrel {\times 1_\rho} \longrightarrow
( \rho^{r+1},\rho^{r+k+1}) \stackrel {\times 1_\rho } \longrightarrow ...
$$

\noindent having assumed that tensoring on the right with the identity arrow $1_\rho$ is an isometric map. Note that $\oro^0$ is a \sC algebra. The direct sum $^0\oro := \oplus_k \oro^k$, with product given by composition of arrows and involution induced by the analogue categorical structure, is a {\it *}-algebra. On $^0\oro$ is naturally defined a circle action, assigning to each $z$ in ${\bT}$ the map 

\begin{equation}
\label{1}
\varphi_z(t_k) := z^k t_k 
\end{equation}

\noindent where $t_k \in \oro^k$. Such a structure is called a ${\bZ}$-graded \sC algebra. It can be shown that there exists a unique \sC norm on $^0\oro$ extending the \sC norm on $\oro^0$ such that (\ref{1}) extends to an automorphic action; the completition $\oro$ is called the DR-algebra associated to $\rho$. The isomorphism class of $\oro$ does not change as $\rho$ varies in its unitary equivalence class (we say that $\rho$ is unitarily equivalent to $\sigma$ if there exists $u \in (\rho,\sigma)$ such that $u^* \circ u = 1_\rho ,u \circ u^* = 1_\sigma$): in fact, for each suitable $r,k\in {\bZ}$ we have the Banach space isomorphisms 

\begin{equation}
\label{2}
\begin{array}{cccc}
\varphi_u^{r,k}: & (\rho^r,\rho^{r+k}) & \ra & (
\sigma ^r,\sigma ^{r+k}) \\  
& t & \mapsto & u^{\otimes^{r+k}} t {u^*}^{\otimes^r} 
\end{array}
\end{equation}

\noindent (in the following we will adopt this notation also when, more generally, 
$u \circ u^* \neq 1_\sigma$) realizing an isomorphism of ${\bZ}$-graded \sC algebras and hence a grading preserving isomorphism at the level of the associated \sC algebras. In particular for $\rho =\sigma $ the algebra $\oro$ is acted upon by the group ${\bf U}_\rho $ of unitary arrows in $(\rho,\rho)$. $\oro$ is also naturally equipped with a {\em canonical endomorphism} $\widehat{\rho}$ obtained by tensoring on the left with the identity arrow $1_\rho $. Note that $\widehat{\rho}$ commutes with the action of ${\bf U}_\rho$.

A similar construction can be done in the more general setting of strict semitensor \sC categories (see \cite {DPZ97}: roughly speaking, at the level of arrows, only the operation of tensoring on the right with an identity arrow is allowed), for example the category of finitely generated Hilbert bimodules.

Examples of DR-algebras are the Cuntz algebras $\mO_n$, associated to
$n$-dimensional Hilbert spaces with $n < \infty$ (at infinite dimension the Cuntz algebra $\mO_\infty $ is strictly contained in the DR-algebra associated to the Hilbert space, see \cite{CDPR94}), the Pimsner \sC algebras \cite{Pim97,DPZ97} associated to Hilbert bimodules and the algebras $\mO_G$ associated to representations of a compact Lie group $G$ \cite{DR87}.

\subsection{Continuous Fields of \sC algebras with a circle action.}

The following definition for continuous fields of Banach spaces and \sC algebras appears in \cite {KW95}. Although less general, it has the advantage of being less technical in comparison with the original one \cite{Dix,DD63}. The two definitions coincide for a compact base space, that is the case in which we are interested.

\begin{defn}
A continuous field of Banach spaces over a locally compact Haussdorff space $\Omega$ is a pair $( X ,\pi_\om : X \ra X_\om )$, where $X$ is a Banach space and $(\pi_\om)_\om$ is a family of surjective Banach space maps indexed by the points of $\Omega $. The Banach spaces $X_\om$ are called the fibers of $X$ in $\om $. For each $x$\ in $X$ the family $( x_\om := \pi_\om (x))_\om$ is called the vector field associated to $x$. We require the following properties:

\begin{itemize}
\item [(1) \quad]  the norm function $\om \mapsto \left\| x_\om \right\|$ belongs to $C_0(\Omega)$ for each $x \in X$

\item [(2) \quad]  for each $x$\ in $X$, $\left\| x\right\| =\sup \limits_\om \left\| x_\om \right\| $

\item [(3) \quad]  $X$ is a $C_0(\Omega)$-module w.r.t. pointwise multiplication by continuous, vanishing at infinity maps over $\Omega $.

\end{itemize}
\end{defn}

The corresponding definition for continuous fields of \sC algebras is analogous. If $(\mA, \mA \ra \mA_\om)$ is a continuous field of \sC algebras then $\mA$ is called the {\em C*-algebra of the field}. For further basic notions and terminology (morphisms, local triviality) see the references above.

We now give a characterization of locally trivial continuous fields of \sC algebras with constant fiber $\mA_\om \equiv \mA$ in terms of the cohomology set $H^1( \Omega ,aut\mA)$ (see \cite{Kar} for the definition), where the group of automorphisms $aut\mA$ is endowed with the pointwise convergence topology. These classical ideas were first stated in \cite {DD63}, \S 26 for continuous fields having as fiber the compact operators, but the same argument works in the general case. The proof of the next theorem, that we omit, is based on usual arguments of transition functions for fiber bundles and makes use of a clutching lemma (\cite {Dix} 10.1.13). The hypothesis of connectedness is necessary to fix the fiber, and the compactness of $\Omega$ can be replaced with paracompatness using the definition by Dixmier and Douady.

\begin{thm}
\label{th21}
Let $\Omega$ be a compact connected Haussdorff space and $\mA$\ a C*-algebra. Then there is a bijective correspondence between the set of isomorphism classes of locally trivial continuous fields of \sC algebras with fiber $\mA$ over $\Omega$\ and the cohomology set $H^1(\Omega ,aut\mA)$. 
\end{thm}

By a general result (\cite{Ste}) we have $H^1(S^n,G) =\pi_{n-1} (G) / \pi_0 (G)$ for each topological group $G$, where $S^n$ denotes the $n$-sphere and $\pi_0 (G)$ is the set of connected components of $G$. Thus the classification of locally trivial continuous fields of \sC algebras over the $n$-spheres is equivalent to the homotopy theory of the group of automorphisms of the fiber. There are several results in this direction in the case of AF-algebras \cite {Nis88, Tho87, Nis92}\footnote{The author thanks P. Goldstein for drawing his attention to these references.}. The next corollaries will be used and illustrated by explicit examples in the sequel.

\begin{cor}
\label{cor21}
With the notation of the previous theorem, let $(\mA,\alpha)$ be a \sC dynamical system over a topological group $G$. Then $\alpha$ induces a map from $H^1(\Omega,G)$ into the set of isomorphism classes of locally trivial continuous fields of \sC algebras with fiber $\mA$.
\end{cor}

\noindent By the previous corollary, the existence of the automorphic action (\ref{1}) and the well known isomorphism $H^1(\Omega,\bT) \simeq H^2(\Omega,\bZ)$ we deduce the following result:

\begin{cor}
\label{cor31}
Let $\mA$ be a \sC algebra carrying a strongly continuous ${\bT}$-action by automorphisms (as in particular a DR-algebra). Then each element of the Cech cohomology group $H^2(\Omega,{\bZ})$ defines a locally trivial continuous field over $\Omega$ with fiber $\mA$.
\end{cor}

\section{Extensions of \sC Categories.}

Let $\mC$ be a \sC category and $\mA$ a \sC algebra with identity $1$. Given a pair $\rho,\sigma$ of objects in $\mC$ we consider the space of arrows $(\rho,\sigma)$; it has the following structure of a Hilbert $(\sigma ,\sigma)$--$(\rho ,\rho)$--bimodule:

\begin{equation}
\left\{ \begin{array}{l}l,s\mapsto l \circ s \\ s,m\mapsto s\circ m \\ 
\left\langle s,t\right\rangle _R:=s^{*}\circ t\in ( \rho ,\rho ) \\
\left\langle s,t\right\rangle _L:=s\circ t^{*}\in ( \sigma ,\sigma) 
\end{array}
\right.
\end{equation}

\noindent where $s,t \in (\rho,\sigma)$, $m \in (\rho,\rho)$, $l \in (\sigma,\sigma)$. As for example in \cite{Lan} we construct the exterior tensor product $\mA \otimes (\rho,\sigma)$, that is a right Hilbert \sC module over the spatial tensor product $\mA \otimes (\rho,\rho)$. In order for more concise notation we write $(\rho,\sigma)^\mA := \mA \otimes (\rho,\sigma)$, while the symbol $(\rho,\sigma)_a^\mA$ will indicate the algebraic tensor product. We also adopt the Sweedler notation $x: = x_1 \otimes x_2,x_1 \in \mA$, $x_2 \in (\rho,\sigma)$ for elements of $(\rho,\sigma)^\mA$. We now define a composition law

$$
(\sigma,\tau)_a^\mA \times (\rho,\sigma)_a^\mA \ra (\rho,\tau)_a^\mA 
$$

\noindent extending naturally the one defined on $\mC$: 

\begin{equation}
yx := y_1 x_1\otimes (y_2 \circ x_2) \in (\rho,\tau)_a^\mA
\end{equation}

\noindent where $y \in (\sigma,\tau)_a^\mA$, $x \in (\rho,\sigma)_a^\mA$. With an abuse of notation we also define an involution

$$
* : (\rho,\sigma) _a^\mA \ra (\sigma ,\rho) _a^\mA
$$

\noindent by setting

$$
( x_1 \otimes x_2)^* : = x_1^* \otimes x_2^*.
$$

\noindent The following lemma shows that the extended composition and involution are bounded, hence defined on the closures of the algebraic tensor products.

\begin{lem}
Let $x \in (\rho,\sigma)_a^{\mA} , y \in (\sigma,\tau)_a^\mA$. Then

\begin{itemize}
\item [{\em (1)} \quad]  $\left\| x  \right\|^2 = \left\| x^* x \right\| = \left\| xx^* \right\|$
\item [{\em (2)} \quad]  $\left\| yx \right\| \leq \left\| y \right\| \left\| x \right\| $
\item [{\em (3)} \quad]  $\left\| x  \right\| = \left\| x^* \right\| $
\end{itemize}
\end{lem}

\begin{proof}
\

\quad (1) \quad  It is simply the definition of the 2-norm in $(\rho,\sigma)^\mA$ as a Hilbert \sC bimodule.

\quad (2) \quad  It suffices to observe that $x^*(\left\| y^* y \right\| - y^* y)x$ is positive as an element of $(\rho,\rho)^\mA$. 

\quad (3) \quad  The isometry of the involution is obtained by embedding $(\rho,\sigma)^\mA$ {\it via} the composition law into the space of the bounded $(\rho,\rho)^{\mA}$-linear maps of right Hilbert \sC modules from $(\rho,\rho)^\mA$ into $(\rho,\sigma)^\mA$: if $x \in (\rho,\sigma)^\mA$ then $x^*$ is the adjoint of $x$ as Hilbert \sC module map.
\end{proof}

In order to define the extension category we consider the category $\widetilde{\mC}^\mA$ having the same objects of $\mC$ and arrows $(\rho,\sigma)^\mA$. The previous lemma implies that $\widetilde{\mC}^\mA$ is a \sC category; we close for subobjects and obtain in this way a \sC category ${\mC}^\mA$ whose objects are the projections $e \in (\rho,\rho)^{\mA}$. We call ${\mC}^\mA$ the {\em extension of} $\mC$ {\em by} $\mA$. By definition, the spaces of intertwiners in ${\mC}^\mA$ are given by

$$
(e,f) :=\left\{ x \in (\rho,\sigma)^\mA : xe = fx = x \right\}  
$$

\noindent where $e\in (\rho,\rho)^{\mA}, f \in (\sigma,\sigma)^{\mA}$. Also by definition $(e,f)$ is a Hilbert $(e,e)$-$(f,f)$-bimodule. In particular the following bimodule action of the center $\mZ$ of $\mA$ is defined:

$$
\zeta x := ( \zeta \otimes 1_\sigma ) x = x\zeta := x( \zeta\otimes 1_\rho ) , 
$$

\noindent where $\zeta \in \mZ$, $x \in (e,f)$.

In the case in which $\mA$ is non unital we define $\mC^\mA := \mC^{\mA^+}$, where $\mA^+ := \mA \oplus \bC$. This definition is motivated by the fact that spaces of intertwiners of an object of a \sC category have to be {\em unital} \sC algebras.

\begin{prop}
\label{prop_dir_sum}
Let $\mC$ be a \sC category, $\mA$ a unital \sC algebra. Then

\begin{itemize}
\item [{\em (1)} \quad]  ${\mC}^\mA$ has subobjects;
\item [{\em (2)} \quad]  If $\mC$ have direct sums, then so does ${\mC}^\mA$.
\end{itemize}
\end{prop}

\begin{proof}
The first assertion is obvious. For the second one, let $e\in ( \rho ,\rho )^\mA,f \in ( \sigma ,\sigma )^{\mA}$;\ then $1_\tau = v \circ v^* + w \circ w^*$, $1_\rho =v^* \circ v$, $1_\sigma = w^* \circ w$ for some object $\tau$ in $\mC$ and we can define the isometries $v_e := 1 \otimes v \cdot e$, $w_f := 1 \otimes w \cdot f$. Thus the projection $g := v_e v_e^* + w_f w_f^* \in (\tau,\tau)^{\mA}$ is a direct sum of $e$ and $f$.
\end{proof}

Before we state some simple functorial properties of the above construction let us introduce the following terminology: a functor between two \sC categories is said to be a {\em C*-functor} if commutes with the involutions and the maps induced on the spaces of arrows are bounded and linear.

\begin{prop}
\label{prop33}
The operation $\mC , \mA \ra {\mC}^\mA$ of assigning an extension to a \sC category depends covariantly on the \sC algebra and covariantly (contravariantly) on the category, corresponding to covariant or contravariant \sC functors on $\mC$.
\end{prop}

\begin{proof}
Let $\mC$ be a \sC category, $\phi : \mA \ra \mB$ a \sC algebra morphism. We construct a {\em pullback C*-functor} $\phi_{*}:{\mC}^\mA \ra {\mC}^\mB$. Now there exists, for each object $\rho $ in $\mC$, the \sC algebra morphism

$$
\phi \otimes id_{(\rho,\rho)} : (\rho,\rho)^\mA \ra (\rho,\rho)^\mB ;
$$

\noindent we define $\phi_{*}(e) := \phi \otimes id_{(\rho,\rho)}(e)$. Now, as $\left\| \phi \otimes id_{( \rho,\rho ) }( x) \right\| \leq \left\| x\right\| $ for each $x$ in $( \rho ,\rho ) ^{\mA}$, we define the bounded linear maps

$$
\phi_{*}:(\rho,\sigma)^{\mA} \ra (\rho,\sigma)^\mB 
$$

\noindent by setting

$$
\phi_* ( a \otimes s) := \phi (a) \otimes s, 
$$

\noindent By definition $\phi_*(yx) = \phi_*(y) \phi_*(x)$, so that $\phi_*(e,f) \subset (\phi_*(e) , \phi_*(f) )$. Let now $F: {\mC}_1 \ra {\mC}_2$ be a \sC functor; we want to define a \sC functor $F:{\mC}_1^\mA \ra {\mC}_2^\mA$. We observe first that for each object $\rho$ in ${\mC}_1$ the \sC algebra morphism

$$
id_{\mA}\otimes F:( \rho ,\rho ) ^{\mA}\ra( F\rho ,F \rho)^\mB 
$$

\noindent is induced, so we define

$$
F^{\mA}e := (id_{\mA} \otimes F) (e) . 
$$

\noindent Again, as $id_{\mA} \otimes F$ is bounded, we find that the maps

$$
F^{\mA}(x_1 \otimes x_2) := x_1 \otimes F(x_2), 
$$

\noindent where $x := x_1 \otimes x_2 \in (\rho,\sigma)_a^{\mA}$, can be extended to $(\rho,\sigma)^\mA$. It is trivial to check that composition and involution are preserved. The proof for contravariant \sC functors follows the same lines by using \sC algebra antihomomorphisms. 
\end{proof}

\begin{cor}
If $\mC$ has direct sums, every pullback \sC functor $\phi_* : \mC^\mA \ra \mC^\mA$ preserves direct sums in $\mC^\mA$.
\end{cor}

\begin{proof}
Let $e,f$ be objects in $\mC^\mA$. With the notation used in the proof of Prop. \ref{prop_dir_sum}, it suffices to verify that $\phi_*(v_e v_e^* + w_e w_e^*)$ is a direct sum of $\phi_*(e) , \phi_*(f)$. 
\end{proof}

In order to emphasize a topological point of view we will adopt in the following the notation ${\mC}^\Omega$ to indicate the extension of $\mC$ by a commutative \sC algebra $\mA = C(\Omega)$; in the same spirit we write $(\rho,\sigma)^\Omega := (\rho,\sigma)^{C(\Omega)} \simeq C(\Omega,(\rho,\sigma))$ for each pair of objects $\rho ,\sigma $ in $\mC$. Let $\varphi : \Omega' \ra \Omega$ be a continuous map of compact Haussdorff spaces. Applying the previous proposition we get a pullback functor $\varphi_* : {\mC}^\Omega \ra {\mC}^{\Omega'}$. If $\varphi$ is the inclusion $\left\{ \om \right\} \hra \Omega$ we call the associated functor the {\em fiber functor} of ${\mC}^\Omega$ over $\om$ and use the notation $\om_* : {\mC}^\Omega \ra {\mC}^\om$. $\om_*$ associates to objects and arrows of ${\mC}^\Omega $, regarded as continuous maps, their evaluation in $\om$. Note that by definition ${\mC}^\om$ is isomorphic to the closure for subobjects of $\mC$.

Proposition \ref{prop33} fails for multifunctors. Let in fact $F: {\mC}_1 \times {\mC}_1 \ra {\mC}_2$ be a \sC functor; proceeding with the same argument we could try to define a map

$$
F^{\mA} : (\rho,\sigma)_a^\mA \times (\rho',\sigma')_a^\mA \ra ( F(\rho ,\rho ^{\prime }) ,F(\sigma ,\sigma ') )_a^\mA 
$$

\noindent by the expression 

\begin{equation}
\label{4}
F^{\mA}( x_1 \otimes x_2,x_1^{\prime } \otimes
x_2^{\prime}) := x_1 x_1^{\prime } \otimes F( x_2,x_2^{\prime }) .  
\end{equation}

\noindent In order to verify that this map preserves the composition, let $x \in (\rho,\sigma )_a^{\mA}$, $x' \in ( \rho',\sigma')_a^{\mA}$, $y \in (\sigma,\tau)_a^{\mA}$, $y' \in (\sigma',\tau')_a^{\mA}$. Then we have

$$
F^{\mA}( yx,y^{\prime }x^{\prime }) =y_1x_1y_1^{\prime}x_1^{\prime }\otimes
F( y_2,y_2^{\prime }) \circ F(x_2,x_2^{\prime }) ;
$$

\noindent however

$$
F^{\mA}( y,y^{\prime }) F^{\mA}( x,x^{\prime}) =
y_1y_1^{\prime }x_1x_1^{\prime }\otimes F( y_2,y_2^{\prime}) \circ F(x_2,x_2^{\prime })
$$

\noindent so that unless $\mA$ is commutative there is no natural way to define $F^\mA$. We specialize for the moment our discussion to the commutative case, but will return to the general case in the subsequent sections.

\begin{prop}
Let\ $\mC$ be a \sC category, $\Omega$\ a compact Haussdorff space. Then

\begin{itemize}
\item [{\em (1)} \quad]  Each \sC functor $F:{\mC}_1\times {\mC}_1\ra {\mC}_2$ induces a \sC functor $F^\Omega :{\mC}_1^\Omega \times {\mC}_1^\Omega \ra {\mC}_2^\Omega $.

\item [{\em (2)} \quad]  ${\mC}^\Omega $ is a strict tensor \sC category if $\mC$ itself is so; the space of intertwiners of the identity object $\iota_\Omega := 1 \otimes 1_\iota$ in ${\mC}^\Omega $ is $( \iota _\Omega,\iota _\Omega ) \simeq {\mC}( \Omega) \otimes (\iota ,\iota ) $, where $\iota $ is the identity object in $\mC$.

\item [{\em (3)} \quad]  ${\mC}^\Omega $ is symmetric if $\mC$ is symmetric. 

\item [{\em (4)} \quad]  If $\mC$ has conjugates, then so does ${\mC}^\Omega $.

\end{itemize}
\end{prop}

\begin{proof} 
The computation above shows that the map (\ref{4}) preserves the composition of arrows in ${\mC}^\Omega $. Using commutativity of $C(\Omega)$ it is also verified that $F^\Omega $ preserves the involution. Now, thanks to the \sC identity, in order to verify that $F^\Omega$ is bounded, it suffices to show this just in the case where the domain object coincides with the codomain. We note that the map $s \mapsto F(s,1_\sigma)$ defines a morphism $(\rho,\rho) \ra (F(\rho,\sigma) , F(\rho ,\sigma))$, that we extend to a morphism $F_\rho^\Omega : (\rho,\rho)^\Omega \ra ( F(\rho,\sigma) , F(\rho,\sigma))^\Omega $ tensoring by the identity on $C( \Omega)$. Thus, for $x \in (\rho,\rho)^\Omega $, the inequality $\left\| F_\rho^\Omega ( x) \right\| \leq \left\| x\right\|$ holds. Now for $x \in (\rho,\rho)_a^\Omega ,y\in ( \sigma ,\sigma )_a^\Omega $ we have

$$
F^\Omega (x,y) =F^\Omega ( x,1_\sigma ) F^\Omega( 1_\rho ,y)
= F_\rho ^\Omega (x) F_\sigma ^\Omega(y) 
$$

\noindent so that

$$
\left\| F^\Omega ( x,y) \right\| \leq \left\| x\right\| \left\|y\right\| . 
$$

\noindent Thus we can extend $F^\Omega$ to $(\rho,\sigma)^\Omega $ and define, for each pair $e,f$ in ${\mC}^\Omega $, the object $F^\Omega (e,f) $ in ${\mC}_2^\Omega $ that can be expressed conveniently by the following formal expression

$$
F^\Omega (e,f) := e_1 f_1\otimes F(e_2,f_2) . 
$$

\noindent We pass now to prove the second point. We define the tensor product on ${\mC}^\Omega $ as the \sC functor

$$
ef := e_1 f_1 \otimes e_2 \times f_2 
$$

\noindent induced by the tensor product on $\mC$; with an abuse of notation, we write

$$
x \times y := x_1 y_1 \otimes x_2 \times y_2 
$$

\noindent to indicate the tensor product of arrows in ${\mC}^\Omega $. By (1), the only properties we have to prove are strict associativity and the existence of an identity object. Strict associativity follows by using the strict associativity of the tensor product of arrows in $\mC$. For the symmetry, we define

$$
\theta_\Omega (e,f) :=1 \otimes \theta (\rho,\sigma) \cdot e \times f
$$

\noindent (in the previous expression we regarded $e,f$ as the identity arrows of $(e,e),(f,f)$); the properties of symmetry are verified with direct computations. For the existence of conjugates, we observe that $\widetilde{\mC}^\Omega $ as defined at the beginning of this section obviously has conjugates. Then we apply \cite{LR97}, Th. 2.4, where it is shown that each subobject of an object having conjugates itself has conjugates. 
\end{proof}

The following result emphasizes a geometrical interpretation of ${\mC}^\Omega$. Roughly speaking we can say that while a \sC category is the categorical analogue of a \sC algebra, ${\mC}^\Omega$ is an example of the categorical analogue of a continuous field of \sC algebras.

\begin{prop}
\label{prop35}
Let $\mC$ be a \sC category closed for subobjects, $\Omega $ a compact Haussdorff space. Then

\begin{itemize}
\item [{\em (1)} \quad]  Each object $e$ in ${\mC}^\Omega $ defines, via the fiber functors $\om
_{*}:{\mC}^\Omega \ra \mC$, a family $\left\{ e_\om \right\}_{\om \in \Omega }$ of objects in ${\mC}$, called the fibers of $e$;

\item [{\em (2)} \quad]  The involution on ${\mC}^\Omega $ is defined fiberwise, i.e. $* \circ \om_* = \om_* \circ * $ for each $\om \in \Omega$;

\item [{\em (3)} \quad]  For each $e,f$ objects in ${\mC}^\Omega $, the space of arrows $(e,f)$
defines a locally trivial continuous field of Banach spaces 

\begin{equation}
\label{5}\left\{ (e,f) ,\om_{*} : (e,f) \ra (e_\om,f_\om) \right\} ,
\end{equation}

\noindent over $\Omega$, with fibers the spaces of intertwiners $( e_\om , f_\om)$ of the fibers of $e,f$ in $\mC$.

\item [{\em (4)} \quad]  If $\mC$ is strict tensor the tensor product defined on ${\mC}^\Omega $ induces on the spaces of arrows morphisms of continuous fields of Banach spaces.
\end{itemize}
\end{prop}

\begin{proof} 
The first and second points are obvious. For the third point, note that $e \in (\rho,\rho)^\Omega , f \in (\sigma,\sigma)^\Omega$ can be regarded as continuous functions $\om \mapsto e(\om) ,f(\om)$. Then we consider in $\mC$ the subobjects $e_\om := \om_*(e) , f_\om := \om_*(f)$ defined by $e(\om) , f(\om)$, and prove that (\ref{5}) is a continuous field of Banach spaces. If $x \in (e,f) \subset (\rho,\sigma)^\Omega$ then the continuous function $\om \mapsto x(\om)$ is defined; we put $x_\om := \om_*(x) \in (e_\om,f_\om)$ and observe that $(e_\om,f_\om)$ is isometrically isomorphic as a Banach space to $(e(\om),f(\om)) \subset (\rho ,\sigma)$ {\em via} the map $s \mapsto w_\om s v_\om^*$, where $s \in (e_\om,f_\om) , v_\om^* v_\om = 1_{e_\om}, v_\om v_\om^* = e_\om , w_\om^* w_\om = 1_{f_\om } , w_\om w_\om^* = f_\om$. Thus 

\begin{equation}
\label{6}
\left\| x_\om \right\| = \left\| x(\om) \right\| 
\end{equation}

\noindent and the function $\om \mapsto \left\| x_\om \right\| $ is continuous. The equality (\ref{6}) implies also that $\left\| x\right\| =\sup \limits_\om \left\|x_\om \right\| $. We already know that $(e,f)$ is a Banach $C(\Omega)$-module, so have proved that
(\ref{5}) is a continuous field of Banach spaces. We postpone the proof of local triviality until the next lemma. Finally, we have to prove that the map $x,y \mapsto x \times y$, $x \in ( e,e')$, $y \in (f,f')$ preserves the fibers, i.e. $\om_*( x \times y) = \om_*(x) \times \om_*(y)$, but this is obvious regarding $x,y$ as continuous map and recalling the definition of the fiber functor.
\end{proof}

\begin{lem}[Local Triviality]
Let $e,f$ be objects in ${\mC}^\Omega$. Then

\begin{itemize}
\item [{\em (1)} \quad]  If the maps $\om \mapsto e(\om ),f(\om)$ are constant there exist objects $\rho_0,\sigma_0$ in $\mC$ such that there is an isomorphism of continuous fields of Banach spaces $(e,f) \simeq (\rho_0,\sigma_0)^\Omega$;

\item [{\em (2)} \quad]  Each point $\om$ in $\Omega$ admits a neighborhood $U$ such that there is an isomorphism of continuous fields of Banach spaces $\left. (e,f) \right|_U \simeq
(e_\om,f_\om)^U$.

\end{itemize}
\end{lem}

\begin{proof}
We define $\rho_0 := e_\om , \sigma_0 := f_\om$ . By hypothesis there exist isometries $v \in (\rho,\rho_0) , w \in (\sigma,\sigma_0)$ intertwining $e,f$ with the units arrows $1_{\rho_0},1_{\sigma_0}$. Thus the map $x \mapsto w x v^{*} , x \in (e,f)$, defines the isomorphism between $(e,f)$ and $(\rho_0,\sigma_0)^\Omega$. For the second point, since $e \in (\rho,\rho)^\Omega ,f \in (\sigma,\sigma)^\Omega$, for each $\om$ there exists a closed neighborhood $U$ such that $\left\| e(\om) - e(\om') \right\|$, $\left\| f(\om) - f(\om') \right\| \leq \delta < 1 , \om' \in U$. This implies that $\left. e \right|_U , \left. f \right|_U$ and the constant maps $e_U := e(\om) , f_U := f(\om)$ are Murray-Von Neumann equivalent as elements of the \sC algebras $(\rho,\rho)^U , (\sigma,\sigma)^U$. Applying (1), since $e_U,f_U$ are constant maps, there must be an isomorphism $\left. (e,f) \right|_U \simeq (e_\om,f_\om)^U$.
\end{proof}

Some examples of extension categories follow.

\begin{ex} The category of finite dimensional Hilbert spaces is a tensor \sC category if endowed with the usual tensor product. To get an equivalent {\em strict} tensor \sC category we proceed as in \cite {DR89} and consider the category ${\bf Hilb}$ with objects the positive integers $n \in \bN$ and arrows the spaces of $n\times m$ matrices $(n,m) := {\bM}_{n,m}$. The tensor product on objects is the multiplication, while on arrows it is defined by the "lexicographical" product for matrices (see for example \cite {DR89} \S 1.4). The extension ${\bf Hilb}^\Omega $ is the category of projective, finitely generated $C(\Omega)$-modules hence, by the Swan-Serre theorem \cite {Swa62}, equivalent to the category of Hermitian vector bundles over $\Omega $. In the non commutative case we get the finitely generated projective (right) Hilbert $\mA$-modules.
\end{ex}

\begin{ex} Let $\widehat{G}$ be the category of unitary, finite dimensional, continuous representations of a compact group $G$. Each object $(n,u) $ in $\widehat{G}$ is a topological group map $G \ra {\bf U}_n$ into the unitary group of the $n$-dimensional matrices. An object in $\widehat{G}^\Omega$ is a projection $e_u$ in ${\bM}_n \otimes C(\Omega)$ that commutes with the action of $G$ induced by $u$ and, given $(n,u)$, $(m,v)$ objects in $\widehat{G}$, the intertwiners $t\in (e_u,f_v) $ are continuous functions from $\Omega $ into the space of $n\times m$ matrices and satisfy the relations $t = t e_u = f_v t , v_g t = t u_g$ for each $g \in G$. Hence we get the category of $G$-Hermitian vector bundles with trivial action of $G$ over the base space $\Omega$ (\cite{Ati} \S 1.6, \cite{DR89} \S 1.6). The non commutative analogue is given by the finitely generated projective Hilbert $\mA$-modules carrying a unitary representation of $G$.
\end{ex}

\begin{ex} We can generalize the previous example by considering the category of finite dimensional
(co)representations of a unital Hopf \sC algebra $(\mA,\Delta)$, i.e. a compact quantum group in the sense of Woronowicz (see for example the introductory paper \cite{MVD99} and related references). This is a strict tensor \sC category having direct sums, subobjects and conjugates. The objects are given by pairs $(n,u)$, where $u \in {\bM}_n \otimes \mA$ satisfies

$$
\Delta (u_{lm}) = \sum \limits_k u_{lk} \otimes u_{km} , 
$$

\noindent while an arrow $x$ from $(n,u)$ to $(m,v)$ is an $n \times m$ matrix such that $v( 1 \otimes x) = (1 \otimes x) u$. An object in the extension category over $\Omega$ is then a Hermitian vector bundle $\mE$ carrying an injective comodule action by $\mA$, in the sense explained below: the tensor product $\mE \otimes \mA$ is naturally endowed with the structure of a vector $\mA$-bundle in the sense of \cite{MF80}. Then a vector bundle morphism $\varphi_u : \mE \ra \mE \otimes \mA$ is induced, defined fiberwise by setting

\begin{equation}
\label{6.1}
\varphi_u(\xi_\om) := u(\xi_\om \otimes 1) ,
\end{equation}

\noindent where $\om \in \Omega , \xi_\om \in \mE_\om \simeq \bC^n$. Now tensoring again with $\mA$ we obtain a vector $\mA \otimes \mA$-bundle and $id_{\mE} \otimes \Delta$ induces a vector bundle morphism $\mE \otimes \mA \ra \mE \otimes \mA \otimes \mA$. Then we have the identity

$$
(id_{\mE} \otimes \Delta ) \varphi_u = (\varphi_u \otimes id_{\mA}) \varphi_u .
$$

\noindent Note that similar comodule actions over vector bundles can be defined more generally for continuous fields of Hopf \sC algebras.
\end{ex}

\begin{ex} \label{enda} Let $\mA$ be a unital \sC algebra with trivial center. Then the category ${\bf end}\mA$ \cite{DR89}, having as objects the endomorphisms of $\mA$ and arrows 

$$
(\rho,\sigma) = \left\{ t \in \mA : t \rho(a) = \sigma (a) t ,a \in \mA \right\}, 
$$

\noindent is a strict tensor \sC category if endowed with the tensor structure given by the composition of endomorphisms and tensor product on arrows defined by $r\times s:=r\rho (s) = \rho'(s) r$, $r \in (\rho,\rho') , s \in (\sigma,\sigma')$. We consider the extension category ${\bf end}^\Omega \mA$ having as objects the projections in $(\rho,\rho)^\Omega $. Let $\rho \in {\bf end}\mA$ such that the space of intertwiners $(\iota,\rho)$ is a not null Hilbert space in $\mA$ (in the sense of \cite{DR88}); given an object $e$ in ${\bf end}^\Omega \mA$ we can regard it as a projection on $(\iota,\rho)^\Omega $ and get in this way the continuous field of Hilbert spaces ${\mathcal H}_e:=e(\iota ,\rho)^\Omega $, embedded as a $C(\Omega)$-bimodule in $C(\Omega ) \otimes \mA$ in the sense of \cite {DPZ97}. Now if $(\iota ,\rho) $ is finite dimensional ${\mathcal H}_e$ is finitely generated as $C(\Omega)$-bimodule; thus, given a set of generators $\left\{ \psi _k\right\} $ of ${\mathcal H}_e$, there exists an endomorphism $\rho_e$ of $C(\Omega) \otimes \mA$ defined (independently from the choice of the $\psi _k$'s) by the map $a\mapsto \sum\limits_k\psi _ka\psi _k^{*},a\in C( \Omega ) \otimes \mA$. Note that $\rho_e$ is a $C(\Omega)$-module map.
\end{ex}

\bigskip

Our purpose is now to give a description of the objects of ${\mC}^\Omega$ in terms of principal bundles, giving a categorical version of well known facts about vector bundles, such as the Swan-Serre Theorem. In the following, $\mC$ will denote a \sC category closed for subobjects. In order to simplify the exposition we will assume that $\Omega$ is connected.

\begin{lem}
\label{lem_fib}
Let $e$ be an object in ${\mC}^\Omega$. Then for each $\om,\om' \in \Omega$ the associated fibers $e_\om , e_{\om'}$ are unitarily equivalent in $\mC$.
\end{lem}

\begin{proof}
By compactness of $\Omega$ and continuity of $e$ there exists, using the argument of the local triviality lemma, a finite open cover $\left\{ \Omega _i\right\} $ such that $e(\om)$ and $e(\om')$ are Murray-Von Neumann equivalent as elements of $(\rho,\rho)$ if $\om$ and $\om'$ belong to the same $\Omega_i$. Now there must exist some $\om''$ belonging to $\Omega_{ij}$, so that the Murray-Von Neumann equivalence class of $e(\om)$ is constant as $\om$ varies in $\Omega_i \cup \Omega_j$. $\Omega$ being connected, we obtain that $e(\om)$ and $e(\om')$ are equivalent for each $\om,\om'$, so there exists an isometry $v_{\om',\om}$ such that $v_{\om',\om }v_{\om',\om}^{*} = e(\om') , v_{\om',\om}^{*} v_{\om',\om} = e(\om)$. The claim of the lemma now follows observing that, by definition of $e_\om $ , $e(\om) = w_\om^{*} w_\om ,1_{e_\om} = w_\om w_\om^{*}$, so that $u_{\om',\om } := w_{\om'} v_{\om',\om} w_\om^*$ is a unitary arrow in $(e_\om , e_{\om'})$. 
\end{proof}

\begin{defn}
Let $e \in (\rho,\rho)^\Omega$ be an object in ${\mC}^\Omega$. A trivialization of $e$ is a family $(\rho_0,\left\{\Omega_i,v_i\right\}_i)$ where $\rho _0$ is an object in $\mC$, $\left\{ \Omega _i\right\} $ is an open cover of $\Omega $ and $v_i\in( \rho ,\rho _0) ^{\overline{\Omega _i}}$ satisfy the relations $\left. e\right|_{\Omega _i}=v_i^{*}v_i,1_{\rho _0} = v_i v_i^{*}$.
\end{defn}

\begin{lem}
Let $e$ be an object in ${\mC}^\Omega $.

\begin{itemize}
\item [{\em (1)} \quad]  $e$ admits a trivialization $( \rho _0,\left\{ \Omega_i,v_i\right\} _i) $.

\item [{\em (2)} \quad]  There is a continuous map $\alpha_{ij} : \Omega_{ij} := \Omega_i \cap \Omega _j \ra ( \rho_0,\rho_0)$ such that $\alpha _{ij} \alpha_{ij}^{*} = \alpha_{ij}^{*} \alpha _{ij} = 1_{\rho_0}$ in $\Omega_{ij}$, $\alpha_{ij} \alpha_{jk} = \alpha_{ik}$ in $\Omega _{ijk}$.

\item [{\em (3)} \quad]  The equivalence class of $( \left\{ \Omega _i\right\} ,\left\{\alpha _{ij}\right\})$ in $H^1( \Omega ,{\bf U}_{\rho_0}) $\ does not depend on the choice of the
trivialization of $e$.

\end{itemize}
\end{lem}

\begin{proof}
The existence of the trivialization follows from local triviality and the previous lemma. Now given the trivialization $(\rho_0,\left\{ \Omega_i,v_i\right\}_i)$, $\alpha_{ij} :=v_iv_j^{*}$ then satisfies the conditions required in the second point. Regarding the last point, we fix a fiber $\rho _0$ of $e$ and consider two trivializations $(\rho_0,\left\{ \Omega_i,v_i\right\}_i)$, $(\rho_0,\left\{ \Omega_l^{\prime},v_l^{\prime}\right\}_l)$, defining $\alpha_{ij} := v_i v_j^{*}$, $\alpha_{lm}^{\prime} := v_l^{\prime} v_m^{\prime *}$. Then $\alpha_{ij} = u_{il} \alpha_{lm}^{\prime} u_{jm}^{*}$, where $u_{il} = v_iv_l^{\prime *}$. 
\end{proof}

\begin{prop}
\label{thm_co}
Let $\mC$ be a \sC category closed for subobjects and $\Omega$ a connected, compact Haussdorff space. Then for each object $e$ in ${\mC}^\Omega$ there exists an object $\rho_0$ in $\mC$, defined up to unitary equivalence, and an element $\delta(e)$ in the cohomology set $H^1(\Omega,{\bf U}_{\rho_0})$. Two objects $e$ and $f$ in ${\mC}^\Omega$ are unitarily equivalent if and only if $\delta(e) = \delta(f)$.
\end{prop}

\begin{proof}
Define $\delta(e)$ to be the class in $H^1(\Omega ,{\bf U}_{\rho _0}) $ of the cocycles coming from trivializations of $e$. If $e=vv^*$ and $f=v^*v$ are unitarily equivalent then we can choose the same fiber $\rho_0$ for $e$ and $f$, and if $(\left\{ \Omega _i\right\} ,\left\{ \alpha _{ij}:=v_iv_j^*\right\} )$ is a cocycle for $e$, then $( \rho _0,\left\{ \Omega _i,v_iv\right\}_i) $ is a trivialization of $f$ defining the cocycle $ v_ivv^{*}v_j^{*}=\alpha _{ij}$ so that {\it $\delta ( e) =\delta (f) $}. Let now $e$ and $f$ be such that {\it $\delta (e) =\delta ( f) $.} Then there are two equivalent cocycles $ ( \left\{ \Omega _i\right\} ,\left\{ \alpha _{ij}:=v_iv_j^{*}\right\} )$, $( \left\{ \Omega _l^{\prime }\right\} ,\left\{ \alpha_{lm}^{\prime }:=v_l^{\prime }v_m^{\prime *}\right\} )$ associated respectively to $e$ and $f$, and there exist continuous maps $u_{il}$ from $ \Omega _i\cap \Omega _l^{\prime }$ into ${\bf U}_{\rho _0}$ such that $ v_iv_j^{*}=u_{il}v_l^{\prime }v_m^{\prime *}u_{jm}^{*}$. We define $w_{il} := v_l^{\prime ~*}u_{il}^{*}v_i\in \left. ( e,f) \right|_{\Omega _i\cap \Omega _l^{\prime }}$; then $w_{il}^{*}w_{il}=\left. e\right| _{\Omega _i\cap \Omega _l^{\prime
}},w_{il}w_{il}^{*}=\left. f\right| _{\Omega _i\cap \Omega _l^{\prime }}$ and $\left.
w_{il}\right|_{\Omega_{ij} \cap \Omega_{lm}^{\prime}} = \left. w_{jm}\right|_{\Omega_{ij} \cap \Omega_{lm}^{\prime}}$ so that we can clutch together the $w_{il}$'s and obtain an isometry $w$ intertwining $e$ and $f$. 
\end{proof}

\section{Continuous Fields of DR-Algebras.}

In this section we prove the main result about the interpretation of DR-algebras associated to objects in ${\mC}^\Omega$ as continuous fields of \sC algebras. We state the result for ${\mC}^\Omega $, but the argument works without substantial changes for 'continuous fields of \sC categories' satisfying the properties stated in Proposition \ref{prop35}. We will assume in the rest of this section that $\mC$ has subobjects. The following lemma states a sufficient condition to be the tensoring with the unit arrow an isometric map; recall that we need this property to construct the inductive limits structures in the DR-algebra.

\begin{lem}
Let $e \in (\rho,\rho )^\Omega$ be an object in ${\mC}^\Omega$. If $\left\| s \right\| =\left\| s \times 1_{e_\om} \right\|$ for each fiber $e_\om$ of $e$, $s \in (\sigma,\sigma')$, then $\left\| x \times e \right\| = \left\| x \right\|$ for each $x \in (f,f') \subset (\sigma,\sigma')
^\Omega$.
\end{lem}

\begin{proof} 
By Proposition \ref{prop35} follows that for each $\om \in \Omega $ there exist isometries $v_\om ,w_\om $ providing an (isometric) \sC algebra isomorphism

$$
\begin{array}{c}
( (fe) (\om ),( f^{\prime }e) (\om)) \ra ( f_\om e_\om ,f_\om^{\prime }e_\om ) \\ 
y(\om )\mapsto v_\om \times w_\om \cdot y( \om ) \cdot v_\om ^{*}\times w_\om ^{*}
\end{array}
$$

\noindent so that

$$
\begin{array}{ll}
\left\| x(\om)\times e(\om)\right\| ^2 & =\left\| x(\om
)^{*}x(\om )\times e(\om )\right\| \\  
& =\left\| v_\om \times w_\om \cdot x(\om )^{*}x(\om )\times
e(\om )\cdot v_\om ^{*}\times w_\om ^{*}\right\| \\
& =\left\| v_\om x(\om )^{*}x(\om )v_\om ^{*}\times 1_{\rho
_\om }\right\| \\
& =\left\| v_\om x(\om )^{*}x(\om )v_\om ^{*}\right\| =\left\|
x(\om )\right\| ^2. 
\end{array}
$$

\noindent This proves the lemma. 
\end{proof}

In the case in which the isometry of tensoring on the right with a unit arrow is not verified we can quotient as in \cite {DR89}, \S 4 for each fiber, so that the hypothesis of the previous lemma is satisfied. In the following we assume that $e$ satisfies the hypothesis of the previous lemma. We start by stating a triviality result.

\begin{lem}
If $e$ is a constant map then there is a \sC algebra isomorphism 
$\oex \cong C(\Omega) \otimes \mO_{\rho_0}$.
\end{lem}

\begin{proof}
By the first point of the local triviality lemma there exists an isometry $v \in (\rho,\rho_0)$ and a family of Banach space isomorphisms $\varphi_v^k := (\varphi_v^{r,k})_r$, defined by

$$
\varphi_v^{r,k}( x) :=v^{r+k}xv^{*r}\in ( \rho _0^r,\rho_0^{r+k}) 
$$

\noindent where $x \in (e^r,e^{r+k})$. Hence $\varphi_v^k$ defines an isomorphism of ${\bZ}$-graded \sC algebras 

$$
^0\oex \cong \bigoplus \limits_k( C(\Omega) \otimes \mO_{\rho_0}^k). 
$$
\end{proof}

\begin{lem}
For each $k\in $\ ${\bZ}$,\ $\oex^k$\ has the structure of a locally trivial continuous field of Banach spaces $( \oex^k,\om_{*} : \mO_e^k \ra \mO_{e_\om}^k)$.  
\end{lem}

\begin{proof}
Tensoring on the right with {\it $e$} gives the isometric embedding of continuous fields of Banach spaces

$$
\begin{array}{ccc}
( e^r,e^{r+k}) & \ra & ( e^{r+1},e^{r+k+1}) \\ x & \mapsto & x\times
e  \end{array}
$$

\noindent and we have the commutative diagram

\begin{center}
$\xymatrix{
           ( e^r,e^{r+k} )
		   \ar[d]_{\times e}
		   \ar[r]^{\om_{*}}
		 & ( e_\om^r,e_\om^{r+k} )
		   \ar[d]^{ \times 1_{\rho ( \om ) } } \\
		   ( e^{r+1},e^{r+k+1} )
		   \ar[r]^{\om_{*}}
		 & ( e_\om ^{r+1} , e_\om^{r+k+1} ) 
}$ \\
\end{center}

\noindent so that the evaluation over $\om$ is well defined for $\oex^k$:

$$
\om_{*} : \oex^k \ra \mO_{e_\om}^k .
$$

\noindent The properties of continuous fields of Banach spaces are naturally inherited by the continuous fields $(e^r,e^{r+k})$. For local triviality, it suffices to note the compatibility of the local charts defined above with the inductive limit structure. 
\end{proof}

Let now $\Sigma$ be a subset of $\Omega$; we define $^0\mO_{e,\Sigma} := \left\{ \mO_{e,\Sigma }^k\right\} _{k\in {\bZ}}$, where

$$
\mO_{e,\Sigma }^k := \left\{ (x_\om)_{\om \in \Sigma} \in \prod \limits_{\om \in \Sigma } \mO_{e_\om}^k , x \in \mO_e^k \right\} . 
$$

\noindent Since $\left\| x_\om \right\| \leq \left\| x\right\| $ on each $\mO_{e,\Sigma}^k$ the norm $\left\| x\right\| _\Sigma := \sup_{\om \in \Sigma} \left\| x_\om \right\| $ is defined, so that $\mO_{e,\Sigma }^k$ has the structure of a ${\bZ}$-graded \sC algebra. We indicate the associated \sC algebra with the notation $\mO_{e,\Sigma }$; it is naturally equipped with the evaluation morphisms $\pi_\om( x) := x_\om \in \mO_{e_\om}$, $x \in \mO_{e,\Sigma}$.

We now observe that the fiber functor $\om_* : {\mC}^\Omega \ra \mC$ induces, by \cite{DR89}, Th. 5.1, a \sC algebra morphism 

$$
\begin{array}{c}
\oex \ra \mO_{e_\om } \\ x \mapsto x_\om 
\end{array}
$$

\noindent By abuse of notation we denote this morphism by $\om_*$. The following lemma will imply the local triviality of the continuous field defined by $\oex$.

\begin{lem}
For each $\om_0$\ in $\Omega$\ there are a neighborhood $U$ of $\om_0$, an isomorphism $\alpha_U : \mO_{e,U} \ra C(\overline{U}) \otimes {\mO}_{e_{\om_0}}$ and a family $(\alpha_\om)_{\om \in U}$ of isomorphisms $\alpha_\om : \mO_{e_\om} \ra \mO_{e_{\om_0}}$ such that, for each $\om \in U$, $\alpha_\om \circ \pi_\om = \om_* \circ \alpha_U$.
\end{lem}

\begin{proof}
The existence of $\alpha_U$ follows by the previous lemmas and the isomorphism is constructed with an isometry $v_U \in (\rho,\rho)^U$ such that $v_U^{*} v_U = \left. e\right|_U$, $v_U v_U^{*} = e(\om_0)$. We have also, for each $\om$, isometries $w_\om$, $w_\om^{\prime}$ implementing the following Banach space isomorphisms

\begin{center}
$\xymatrix{
           ( e_\om^r,e_\om^{r+k} )
		   \ar[d]_{ {\varphi_{w_\om^{*}}^{r,k}} }
		   \ar[r]
		 & ( e_{\om_0}^r , e_{\om_0}^{r+k} )
		\\ ( e^r(\om ) , e^{r+k}(\om )) 
		   \ar[r]^{ \varphi_{v_U( \om ) }^{r,k}}
	     & ( e^r(\om_0) , e^{r+k}(\om_0) )
           \ar[u]_{ \varphi _{w_\om ^{\prime } }^{r,k} }
}$ \\
\end{center}

\noindent the map defined by the upper horizontal arrow providing the isomorphism $\alpha_\om$ that, by construction, satisfies the condition required.
\end{proof}

\begin{cor}
For each $x$ in $\oex$, the norm function $\om \mapsto \left\|x_\om \right\| $\ is continuous, and $\left\| x\right\| =\left\| x\right\| _\Omega :=\sup \limits_\om \left\|x_\om \right\| $. 
\end{cor}

\begin{proof}
Let $\om_0 \in \Omega$; then for $\om$ in a suitable neighborhood of $\om_0$

$$
\left\| x_\om \right\| =\left\| \alpha_\om x_\om \right\| = \left\|(\alpha_Ux) ( \om ) \right\|  
$$

\noindent and the last term is a continuous function as $\om$ varies in $U$. To prove the result about the norm of $x$ it suffices to check that $\left\| \cdot \right\| _\Omega $ has the required property w.r.t. the circle action. But it is obvious that $\left\| \cdot \right\|_\Omega = \left\|\varphi_z\circ \cdot \right\|_\Omega ,z \in {\bT}$.
\end{proof}

The previous results imply the following

\begin{thm}
\label{main-thm}
Let $\mC$ be a strict tensor C*-category closed for subobjects, $\Omega$ a compact
Haussdorff space. If $e$ is an object in $\mC^\Omega$, then $\oex$ is the C*-algebra of a locally trivial continuous field $(\oex,\om_* : \oex \ra \mO_{e_\om})$, where each $\mO_{e_\om}$ is the DR-algebra of an object $e_\om \in \mC$.
If $e$,$f$ are unitarily equivalent objects in $\mC^\Omega$, then the associated C*-algebras are isomorphic as continuous fields. 
\end{thm}

\noindent Recall that by Prop.\ref{thm_co} every object $e$ in $\mC^\Omega$ is described up to unitary equivalence by an element of a cohomology set $H^1(\Omega,{\bf U}_\rho)$. Now (see Sec.2), for each object $\rho$ in $\mC$ the canonical \sC dynamical system $(\varphi,\oro)$ over ${\bf U}_\rho$ is defined. By Cor. \ref{cor21} there is an application $\varphi_{*}$ from $ H^1( \Omega ,{\bf U}_\rho)$ into $H^1(\Omega,aut \oro)$. We have the following:

\begin{prop}
\label{thm_co_alg}
Let $\mC$ be a strict tensor \sC category closed for subobjects, $\Omega$ a connected, compact Haussdorff space.
\begin{itemize} 
\item [{\em (1)} \quad]  The element of $H^1(\Omega,aut \mO_{\rho})$ defined by $\oex$\ is $\varphi_* (\delta(e))$.

\item [{\em (2)} \quad]  For each couple of objects $e,f$ in $\mC^\Omega$, $\oex \cong \mO_f$ iff $\varphi_* (\delta(e)) = \varphi_* (\delta(f))$.

\end{itemize}
\end{prop}

\begin{proof}
Let $(\rho , \left\{ \Omega_i , v_i \right\}_i)$ be a trivialization of $e$ defining the cocycle {\it $\alpha_{ij}$}$:=v_iv_j^{*}$; then the trivialization $\varphi _i:\left. \oex\right| _{\Omega_i}\ra C( \overline{\Omega _i}) \otimes \mO_{\rho}$, with $\varphi_i( x) :=\varphi _{v_i}^{r,k}( x) ,x\in (e^r,e^{r+k}) $, is induced on $\oex$. Thus $( \left\{\Omega _i\right\} ,\left\{ \widehat{\alpha}_{ij}:=\varphi _i\varphi_j^{-1}\right\} ) $ is a cocycle for $\oex$ and for $s\in (\rho^r,\rho^{r+k}) $ we have

$$
\begin{array}{ll}
\widehat{\alpha }_{ij}( \om ) s & = \varphi_{v_i( \om ) }^{r,k}\varphi
_{v_i^{*}( \om ) }^{r,k}(s) \\ 
& =v_i^{r+k}( \om ) v_j^{*~r+k}( \om ) sv_j^r( \om )
v_i^{*~r}( \om )  \\  
& =( v_iv_j^{*}(\om) ) ^{r+k}s( v_iv_j^{*}( \om )
) ^{*~r} \\  
& =\varphi ( \alpha_{ij}(\om) ) s. 
\end{array}
$$

\noindent This completes the proof of the first statement. The second follows by Theorem \ref{th21}.
\end{proof}

The basic example to illustrate the previous results is given by the category of Hermitian vector bundles over a compact Haussdorff space $\Omega$. In this case Thm.\ref{main-thm} implies that the DR-algebra $\mO_{\mE}$ of a Hermitian vector bundle $\mE$ is the \sC algebra of a continuous field of Cuntz algebras over $\Omega$. We give here a brief description of such algebras, postponing a detailed study to a successive work.

We start by considering the spaces of arrows $({\mE}^r,{\mE}^{r+k})$ between tensor powers of $\mE$; recall that by general facts \cite {DR89} 4.3, \cite {LR97} there is a natural identification 

\begin{equation}
\label{7}
( {\mE}^r,{\mE}^{r+k}) \simeq ( \sigma_{\mE}^r,\sigma_{\mE}^{r+k}) :=\left\{ x\in \mO_{\mE} : x \sigma_{\mE}^r(y) = \sigma_{\mE}^{r+k}( y) x,y\in\coe\right\} 
\end{equation}

\noindent where $\sigma_{\mE}$ is the canonical endomorphism on $\mO_{\mE}$. In the case $r=0$ we have the trivial bundle $\iota :=\Omega \times \bC$, so that $( \iota,\mE)$ can be regarded as the module of continuous sections of $\mE$ and is generated by a finite set $\left\{ \psi_l\in ( \iota ,\mE) , l=1,...,N \right\}$. Now, each element of a fiber of $\mE$ being the evaluation of a section at a point $\om$, we have

$$
\xi _\om =\sum\limits_l\psi _l( \om ) \left\langle \psi_l( \om )
,\xi _\om \right\rangle 
$$

\noindent where $\xi _\om \in {\mE}_\om $. This fact is translated in the following relation on the $\psi _l$'s as elements of $\mO_{\mE}$

\begin{equation}
\label{8}
\sum \limits_l \psi_l \psi_l^{*} = 1. 
\end{equation}

\noindent Note that every $(\mE^r,\mE^s)$ has a natural structure of a $C(\Omega)$-$C(\Omega)$-bimodule by considering the coinciding left and right actions given by multiplication with continuous functions over $\Omega$. In particular $(\iota,\mE)$ is a Hilbert $C(\Omega)$-$C(\Omega)$-bimodule. We have the following

\begin{prop}
Let $\mE$ be a Hermitian vector bundle over a compact space $\Omega$, $\left\{ \psi_l \right\}_{l=1}^N$ a finite set of generators for $(\iota,\mE)$. Then: 

\begin{itemize}
\item [{\em (1)} \quad]  Each $({\mE}^r,{\mE}^{r+k})$ is generated as a right $C(\Omega)$-module by the sets 

$$
\left\{ \psi_J \psi _I^{*}, \psi_I \in (\iota,{\mE}^r) , \psi_J \in (\iota,{\mE}^{r+k}) \right\} ; 
$$

\item [{\em (2)} \quad]  $\coe$ is the Pimsner \sC algebra generated by $(\iota,\mE)$ and $C(\Omega)$;

\item [{\em (3)} \quad] $\coe$ is nuclear;

\item [{\em (4)} \quad]  The canonical endomorphism is inner, i.e. for each $x \in \coe$

$$
\sigma_{\mE}(x) =\sum \limits_l \psi _l x \psi _l^{*}. 
$$

\end{itemize}

\end{prop}

\begin{proof}

\
\quad (1) \quad The identity (\ref{8}) implies that $\sum\nolimits_I\psi _I\psi_I^{*}=1$ where $\psi_I := \psi_{i_1} ...\psi_{i_r} \in (\iota,{\mE}^r)$. Hence for $t\in ( {\mE}^r,{\mE}^{r+k}) $ we get $t=\sum\nolimits_{I,J}\psi _I\psi _I^{*}T_{IJ}$, where $T_{IJ}:=\psi_J^{*}t\psi _I \in ( \iota ,\iota ) \simeq C(\Omega)$.

\quad (2) \quad By (1) it follows that $\coe$ is generated as a \sC algebra by $(\iota,\mE)$.

\quad (3) \quad By \cite{DS99} $\coe$ is exact, so that (see \cite{KW95}) for every \sC algebra $\mA$ the maximal and spatial tensor products $\coe \otimes_{max} \mA$, $\coe \otimes \mA$ each have a natural structure of a continuous field with fiber respectively  $\mO_{n(\om)} \otimes_{max} \mA$, $\mO_{n(\om)} \otimes \mA$, where $n(\om)$ is the rank of $\mE$ in $\om$. Since $\mO_{n(\om)}$ is nuclear the maximal and spatial norms coincide over the algebraic tensor product of $\coe$ by $\mA$.

\quad (4) \quad As $(\iota,\mE) = (\iota,\sigma_{\mE})$ we find that $\sigma_{\mE} (t) \psi_l = \psi_l t$ for each $t \in \coe$. Then we apply (\ref{8}).
\end{proof}

We recall that if $\mA \subset \mB$ is an inclusion of \sC algebras a {\em Hilbert $\mA$-module} in $\mB$ \cite {DPZ97} is given by a subspace $\mX$ of $\mB$ closed under right multiplication by elements of $\mA$ and such that $x^*y \in \mA$, $x,y \in \mX$. If there is a finite set $\left\{ x_l\right\} \subset \mX$ such that, with $p=\sum \nolimits_lx_lx_l^{*}$, $x=px$ for each $x\in $ $\mX$ then $p$ is a projection, called the {\em support of }$\mX$, and we say that $\left\{ x_l\right\} $ generates $\mX$ as a Hilbert $\mA$-module. By this definition it follows that $(\iota,\mE)$ is a finitely generated $C(\Omega)$-module in $\coe$ with support the identity. The previous proposition implies

\begin{cor}
Let $\mE,\mF$ be Hermitian vector bundles over $\Omega$.

\begin{itemize}

\item [{\em (1)} \quad]  If $\coe$ is isomorphic to $\mO_\mF$ then $(\iota,\mF)$ appears embedded in $\coe$ as a $C(\Omega)$-module with support the identity.

\item [{\em (2)} \quad]  If $(\iota,\mF) $ is a Hilbert $C(\Omega)$-module in $\mO_\mE$ with support the identity then there is an isometric embedding $\mO_\mF \hra \coe$.

\end{itemize}

\end{cor}

The description of $\coe$ as the \sC algebra generated by the module of continuous sections of $\mE$ supplies a picture in terms of generators and relations. Each set of generators $\left\{ \psi_l\right\}$ of $(\iota,\mE)$, regarded in $\coe$, satisfies the relations

\begin{equation}
\label{9}
\left\{ 
\begin{array}{l} 
\sum \limits_l \psi_l \psi_l^{*} = 1 \\ 
\psi_l^{*} \psi_m=e_{lm}1
\end{array} 
\right. 
\end{equation}

\noindent where $e_{lm} \in C(\Omega)$. Note that (\ref{9}) implies that $(e_{lm}) \in {\bM}_N \otimes C(\Omega)$ is a projection, in fact the projection associated to $\mE$ as a subbundle of $\Omega \times {\bC}^N$. We can use (\ref{9}) as starting point to give an alternative proof of the fact that each projection $e\in {\bM}_N\otimes C( \Omega ) $ defines a continuous field of Cuntz algebras. We proceed on the line of \cite{Lan97} and generalize in the following way: let $\Omega$ be a compact Haussdorff space, ${\mathcal W}_n$ the free *-algebra with identity $1$ generated by $n$ symbols $w_1,...,w_n$. Then we can consider the *-algebra $C( \Omega ) \otimes {\mathcal W}_n$ with certain relations $r_1,...,r_k$ over $1_{C(\Omega)} \otimes w_1,...,1_{C(\Omega ) }\otimes w_n$. Since $r_1,...,r_k$ are relations in $C(\Omega ) \otimes {\mathcal W}_n$, they define, for each $\om \in \Omega$, a set of relations $r_1( \om ) ,...,r_k( \om)$ on ${\mathcal W}_n$, hence a \sC algebra $\mA_\om $obtained by closing ${\mathcal W}_n\left/ \left\{ r_i( \om )\right\} \right. $ with respect a suitable \sC norm $\left\| {}\right\|_\om $ (see \cite {Lan97} for details). Now, by evaluating continuous functions on $\om $, each element $x$ of the *-algebra $( C( \Omega ) \otimes {\mathcal W}_n) \left/ \left\{ r_i\right\} \right. $ defines a vector field $\left\{ x_\om \in {\mathcal W}_n\left/ \left\{ r_i( \om ) \right\} \right. \right\} $. If one shows that the norm function $\om \mapsto \left\| x_\om \right\| _\om $ is continuous, the \sC norm $\left\| x\right\| :=\sup \nolimits_\om \left\| x_\om \right\| _\om $ on $( C( \Omega ) \otimes {\mathcal W}_n) \left/ \left\{ r_i\right\} \right. $ can be defined, obtaining a \sC algebra $\mA$ endowed with natural morphisms $\pi _\om :{\mA\ra \mA}_\om $. Hence $( \mA,\pi _\om :{\mA\ra \mA}_\om ) $ is a continuous field of \sC algebras over $\Omega $. We now regard at (\ref{9}) as a set of relations on $C(\Omega) \otimes {\mathcal W}_N$: then it is easily verified that the relations (\ref{9}) evaluated in $\om $ are equivalent to the relations defining the Cuntz algebra implemented by $n_\om := {\rm rank}(e(\om))$ isometries $\left\{ s_k( \om )\right\} _{k=1}^{n_\om }$ obtained as an orthonormal base of the vector space spanned by the symbols at the point $\om$. Thus the \sC algebra obtained at $\om$ is the Cuntz algebra $\mO_{n_\om}$. The last thing we have to verify to get a continuous field is the continuity of the norm function. Expressing each generator $(\psi_l)_\om$ in terms of the $s_k(\om) $'s we get, using standard calculations, $\left\| (\psi_l)_\om \right\|_\om ^2 = e_{ll}(\om)$, so that the norm function is continuous.

\bigskip

We now give the expression for the cocycle associated to the \sC algebra of a vector bundle $\coe$. As we showed in a more general setting in Prop. \ref{thm_co}, it is well known (cf. for example \cite {Kar}) that the cohomology set $H^1( \Omega , {\bf U}_n)$ classifies the Hermitian vector bundles with rank $n$ over $\Omega $; now if $\mE$ is identified by the ${\bf U}_n$-cocycle $ ( \left\{ \Omega_i\right\} ,\left\{ \alpha_{ij}\right\} ) $ then the $aut\mO_n$-cocycle describing $\coe$ is given by 

\begin{equation}
\label{10}
( \varphi_{*} \alpha_{ij}) (\om) t:= \alpha_{ij}( \om
) ^{\otimes r+k}t \alpha_{ij}(\om)^{\otimes r~*}, 
\end{equation}

\noindent where $t \in (n^r,n^{r+k}) \subset \mO_n$. 

Let us give an explicit example. Recall by Cor. \ref{cor31} that every $z \in H^1(\Omega,\bT)$ defines, for each $n \in \bN$, a locally trivial continuous field $\mF_{z,n}$ having as fiber the Cuntz algebra $\mO_n$. Furthermore it is well known that $z$ can be viewed, up to cocycle equivalence, as a set of transition functions $z_{ij} : \Omega_{ij} \ra {\bT}$ for a line bundle $\mL$ over $\Omega$. Now, applying (\ref{10}) we obtain that $\mF_{z,n}$ is the continuous field associated to $\theta_n \otimes \mL$, where $\theta_n$ is the trivial rank $n$ vector bundle over $\Omega$.

We now prove the existence of non isomorphic vector bundles having isomorphic DR-algebras. By Th. \ref{main-thm} the \sC algebra $\mO_\mL$ associated to a line bundle $\mL$ is a continuous field having as fiber the algebra $C(\bT)$ of continuous functions over the circle, so that is itself commutative. We introduce for $r < 0$ the notation $\mL^r := (\mL^*)^{-r}$. Note that $(\mL,\mL) = C(\bT)$, so that $(\mL^r,\mL^s) = (\iota,\mL^{r-s})$ and

\begin{equation}
\label{lin-alg}
^0\mO_\mL = \bigoplus \limits_{k \in {\bZ}} (\iota,\mL^k) .
\end{equation}

\begin{prop}
Let $\mL$ be a line bundle over a compact Haussdorff space $\Omega$.
\label{prop_iso}
\
\begin{itemize}
\item [{\em (1)} \quad] The spectrum of $\mO_\mL$ is the sphere bundle $\mL_1$ ;
\item [{\em (2)} \quad] $\mO_\mL \simeq \mO_{\mL^*}$
\end{itemize}
\end{prop}

\begin{proof}

\
\quad (1) \quad Let $(\Omega_i,\lambda_{ij}) \in Z^1(\Omega,{\bT})$ be a family of transition functions for $\mL$. Then the sphere bundle $\mL_1$ can be conveniently described by the transition functions $s(\lambda_{ij}) (\om,z) := (\om,\lambda_{ij}(\om)z)$, with $(\om,z) \in \Omega_{ij} \times {\bT}$. Now, note that if $\psi \in (\iota,\mL^k)$ then we have a local description of $\psi$ in terms of continuous functions $\psi_i : \Omega_i \ra \bC$ satisfying the relation 

\begin{equation}
\label{coc_sec}
\psi_i (\om) = \lambda_{ij}^k (\om) \psi_j (\om) .
\end{equation}

\noindent Let now $\pi_i : \left. \mL \right|_{\Omega_i} \ra \Omega_i \times \bC$ be a set of local charts defining the transition functions $(\Omega_i,\lambda_{ij})$. We define continuous functions $\widehat \psi_i : \Omega_i \times \bT \ra \bC$ by setting $\widehat \psi_i (\om,z) := \psi_i (\om) z^{-k}$. Now, as by (\ref{coc_sec}) the identity $\widehat \psi_i \circ s(\lambda_{ij}) = \widehat \psi_j$ holds, we find $\widehat \psi_i \cdot \pi_i = \widehat \psi_j \cdot \pi_j$. We denote by $\widehat \psi$ the continuous function on $\mL_1$ obtained by clutching the maps $\left\{ \widehat \psi_i \cdot \pi_i \right\}$. By extending the map $\psi \mapsto \widehat \psi$ in the obvious way we obtain an isometric $*$-algebra monomorphism $^0\mO_\mL \hra C(\mL_1)$. By applying for every local chart the Stone-Weierstrass theorem on the factor $C(\bT)$ in $C(\Omega_i \times \bT)$ we find that the extended map  $ ~~ \widehat{} : \mO_\mL \ra C(\mL_1)$ is also surjective.

\quad (2) \quad By (\ref{lin-alg}) it follows that $^0\mO_\mL =$ $^0\mO_{\mL^*}$. An alternative way to verify the isomorphism is to observe that complex conjugation on the fibers induces a homeomorphism $\mL_1 \simeq \mL^*_1$: on each local chart we define $\alpha_i (\om,z) := (\om,\overline z)$ and the identity

$$
s(\lambda_{ij}^*) = \alpha_i s(\lambda_{ij}) \alpha_j^{-1}
$$

\noindent holds, so that a global homeomorphism is defined by clutching the maps $\alpha_i$.
\end{proof}

Of course the previous example is very peculiar, our \sC algebra $\mO_\mL$ being abelian. Further examples relative to vector bundles with rank $> 1$, relating isomorphism classes of DR-algebras with topological $K$-theory, can be given over the $n$-spheres. We refer the reader to a work in progress of S. Doplicher, P. Goldstein and the author.

Keeping in mind the basic example of the \sC algebra of a vector bundle it is easy to obtain a similar description for the cases introduced in Sec. 3. If $u : G \ra {\bf U}_n$ is a unitary representation of a compact group $G$ and $\mE$ is the vector bundle defined by $e \in (u,u)^\Omega$ then we obtain an action of $G$ into the unitary arrows in $(e,e)$, so that on each space of intertwiners we have $\varphi_g^{r,k}(t) := g^{\otimes^{r+k}}tg^{* \otimes ^r},t\in (e^r,e^{r+k}) $ (here we write $g$ for $u(g) $). Note that the $G$-action commutes with the canonical endomorphism (see Sec. 2), and we have a grading preserving action by automorphisms of $G$ on $\mO_{\mE}$. The fixed point algebra $\mO_{\mE,u}$ is generated by the intertwiners $(e^r,e^{r+k})_G := \left\{ t \in (e^r,e^{r+k}) : \varphi_g^{r,k}(t) = t \right\}$ and has the structure of a continuous field of \sC algebras with fiber \sC algebras of the type $\mO_G$ defined in \cite{DR87}.

Similarly in the case of a representation $(n,u)$ of a unital Hopf \sC algebra $(\mA,\Delta)$, by extending to $\coe$ the map (\ref{6.1}) we obtain the following coaction (see \cite{Wan99} for a detailed study of similar coactions on Cuntz algebras)

$$
\Phi_u : \coe \ra \coe \otimes \mA .
$$

\noindent As $\coe$ is exact the map above is a morphism of continuous fields of \sC algebras (see \cite{KW95}). We have the fixed point continuous field

$$
\mO_{\mE,u}:=\left\{ t \in \coe : \Phi_u (t) = t \otimes 1 \right\} .  
$$

\section{Non Commutative Extensions.}

In this section we discuss the existence of tensor structures on the extension of a given strict tensor \sC category by a non commutative \sC algebra. In Sec. 3 we showed that in general these structures fail to exist; thus we have to modify our definition for ${\mC}^\mA$, adopting a procedure that is just the abstract version of the left action on a Hilbert \sC bimodule. Constructions of the type we are going to expose have been considered for the case $\mC = {\bf Hilb}$ in \cite{DR89}, \S 1.7 and \cite{SV93,Ruz}, the second ones in the framework of the algebraic quantum field theory.

We now start our general construction. Let $\mC$ be a \sC category, $\mA$ a unital \sC algebra with center $\mZ$. We consider endomorphisms 

\begin{equation}
\label{13}\phi :\mA \ra (\rho,\rho )^{\mA},
\end{equation}

\noindent where $\rho$ is an object in $\mC$, and call them {\em amplimorphisms}. Applying the methods of Sec. 3 we obtain a full subcategory ${\cM}(\mC,\mA)$ of ${\mC}^\mA$, with objects the amplimorphisms $\phi$ and arrows

$$
(\phi,\psi)_{\mM} := (\phi (1),\psi (1)). 
$$

\noindent We also define the not full subcategory ${\cB}(\mC,\mA)$ having the same objects of $ {\cM}(\mC,\mA) $ and arrows

$$
(\phi,\psi) := \left\{ x \in (\phi,\psi)_{\mM} : \psi (
a) x = x \phi (a) ,a\in \mA \right\} . 
$$

\noindent Note that each object $e \in (\rho,\rho)^{\mZ}$ of $ {\mC}^{\mZ}$ defines an amplimorphism $\phi_e(a) := e \cdot a \otimes 1_\rho = a \otimes 1_\rho \cdot e$. Thus there is the following sequence of immersions of \sC categories:

\begin{equation}
\label{14}
{\mC}^{\mZ} \stackrel{\phi_{\cdot}} \hra {\cB}(\mC,\mA) \stackrel{i} \hra {\cM}(\mC,\mA) \hra {\mC}^\mA, 
\end{equation}

\noindent where the third is a full immersion. Each $(\phi,\psi) $ has the structure
of a Banach $\mZ$-bimodule, with coinciding left and right actions
$\psi(a) t = t \phi(a) , t \in (\phi,\psi) , a \in \mZ$. The same definition for $a \in \mA$ provides a structure of Banach $\mA$-$\mA$-bimodule for $(\phi,\psi)_{\mM}$. 

\begin{rem}
Note that we could consider as well endomorphisms of the type $\varphi : \mA \ra (\rho,\rho)^\mB$, obtaining the corresponding structures of Banach $\mA$-$\mB$-bimodules. We denote the corresponding \sC category by $\cB(\mA,\mB;\mC)$. In particular the categorical analogue of a Kasparov module can be defined as a couple $(\varphi,F)$, where $\varphi$ is an object in $\cB(\mA,\mB;\mC)$, $F$ belongs to the ideal of "compact operators" ${\mathcal K}_\varphi := (\iota,\varphi)_\mM (\varphi,\iota)_\mM \subset (\varphi,\varphi)_\mM$ and satisfies the usual properties $[F,\varphi(a)], (F^2 - 1) \varphi(a) , (F-F^*) \varphi(a) \in {\mathcal K}_\varphi$ for every $a \in \mA$.
\end{rem}

\noindent As to be expected, an analogue of Prop. \ref{prop_dir_sum} holds:

\begin{prop}
Let $\mC$ be a \sC category, $\mA$ a unital \sC algebra. Then

\begin{itemize}
\item [{\em (1)} \quad]  ${\cB}(\mA,\mC)$ has subobjects;
\item [{\em (2)} \quad]  If $\mC$ have direct sums, so do ${\cM}(\mA,\mC)$, ${\cB}(\mA,\mC)$.
\end{itemize}
\end{prop}

\begin{proof}

\
\quad (1) \quad  Let $e \in (\phi,\phi)$ be a projection. As $e$ is in the commutant of $\phi(\mA)$ the amplimorphism $a \mapsto \phi(a) \cdot e$ is defined, so that ${\cB}(\mA,\mC)$ has subobjects.

\quad (2) \quad  As ${\cB}(\mA,\mC)$ is a subcategory of ${\cM}(\mA,\mC)$ it suffices to prove the statement only in this first case. Let now $\phi : \mA \ra (\rho,\rho)^\mA$, $\psi : \mA \ra (\sigma,\sigma)^\mA$ be amplimorphism. With the same notations used in the proof of Prop. \ref{prop_dir_sum} we consider a direct sum $\tau$ of $\rho , \sigma$ in $\mC$. Then we define $v_\phi := \phi(1) \cdot 1 \otimes v$, $w_\psi := \psi(1) \cdot 1 \otimes w$ and the map $\delta (a) := v_\phi^* \phi (a) v_\phi + w_\psi^* \psi(a) w_\psi$ from $\mA$ into $(\tau,\tau)^\mA$. We leave to the reader the verification that $\delta$ is an amplimorphism and a direct sum for $\phi$,$\psi$.
\end{proof}

Let now $\mC$ be a strict tensor \sC category with identity object $\iota$. Our aim is to prove that ${\cM}(\mC,\mA)$ can be endowed with a natural strict semitensor structure, while ${\cB}(\mC,\mA)$ becomes a strict tensor \sC category. In order to better illustrate our construction let us give an example: we consider the category $\overline {\bf end}\mA$ having the same objects and arrows as ${\bf end}\mA$ as introduced at the end of Sec. 3, but flipped tensor product $\rho \sigma := \sigma \circ \rho$, $r \times s := \sigma'(r) s = s \sigma(r)$, $r \in (\rho,\rho') , s \in (\sigma,\sigma')$. We embed $\overline {\bf end}\mA$ (not fully) in the category $ {\mM}(\mA)$ having the same objects of $\overline {\bf end} \mA$ and arrows $(\rho,\sigma)_{\mM} := \left\{ x \in \mA : x = \sigma(1) x = x \rho(1) \right\}$. Each $(\rho,\sigma)_{\mM}$ is a $\mA$-$\mA$-Hilbert \sC bimodule with (not commuting) left and right actions given by $\rho$ and $\sigma$. We have also the following semitensor structure:

$$
\left\{ 
\begin{array}{l}
\rho \sigma :=\sigma \circ \rho  \\ 
\Phi_\tau (t) := \tau (t) \in (\rho \tau , \sigma \tau)
_{\mM}, t \in (\rho,\sigma)_{\mM}
\end{array}
\right.  
$$

\noindent By this definition the reason we choose $\overline{\bf end}\mA$ instead of ${\bf end}\mA$ is evident. In ${\bf end}\mA$ the map $\Phi_\tau$ would play the role of a left tensoring with an identity arrow, while we need a right tensoring to get a semitensor structure. The objects in ${\mM}(\mA)$ define naturally Hilbert $\mA$-$\mA$-bimodules $\mX_\rho := (\iota,\rho)_{\mM}$ embedded in $\mA$, on which we can define the tensor product

$$
\mX_\rho \otimes \mY_\sigma :=\left\{ \sigma ( x) y:x\in {\mX}_\rho ,y\in \mY_\sigma \right\} \subset (\iota,\rho \sigma)_{\mM} 
$$

\noindent with $\mA$-valued scalar product $\left\langle x\otimes y,x^{\prime}\otimes y^{\prime }\right\rangle := y^{*}\sigma ( x^{*}x^{\prime }) y^{\prime }$. As expected, we have an embedding of $(\rho,\rho)_{\mM}$ into the \sC algebra of right $\mA$-module maps on $\mX_\rho \otimes \mY_\sigma$ defined by $t \mapsto \sigma (t)$.

$$
$$

We now start the construction of our tensor structures. We have to define a composition law for amplimorphisms; for this purpose we need the universal properties of both the spatial and maximal tensor product, as stated by the following

\begin{lem}
Let $\rho$ be an object of $\mC$ such that the maximal and spatial \sC norms coincide on $(\rho,\rho)_a^{\mA}$. Then for each amplimorphism $\psi : \mA \ra (\sigma,\sigma)^{\mA}$ a morphism $\psi_{\rho,\rho} : (\rho,\rho)^{\mA} \ra (\sigma \rho ,\sigma \rho )^{\mA}$ is defined, given by the expression

\begin{equation}
\label{def_amp}
\psi_{\rho,\rho} ( a \otimes s) := \psi (a)_1\otimes \psi (a)_2 \times s. 
\end{equation}

\end{lem}

\begin{proof}
By property of the spatial tensor product we obtain, tensoring the identity on $\mA$ with the right tensor product with $1_\sigma $, a morphism

$$
i_\rho : (\sigma,\sigma)^\mA \ra (\sigma \rho ,\sigma \rho)^\mA 
$$

\noindent defined by $i_\rho (a\otimes s) :=a\otimes s\times 1_\rho $. By composing with $\psi$ we get a morphism

$$
i_\rho \circ \psi :\mA \ra (\sigma \rho,\sigma \rho)^\mA
$$

\noindent the image of which commutes with that of

$$
j_{\sigma,\rho} : (\rho,\rho) \ra (\sigma \rho,\sigma \rho)^{\mA} 
$$

\noindent defined by $j_{\sigma,\rho}(t) := 1 \otimes 1_\rho \times t$. Thus, by using the universal property of the maximal tensor product, we define $\psi_{\rho,\rho} := ( i_\rho \circ \psi ) \otimes j_{\sigma,\rho}$.
\end{proof}

The hypothesis of the previous lemma is satisfied if, for example, $\mC$ has conjugates (see \cite {LR97}, Lemma 3.2), or $\mA$ is nuclear. In such a case we introduce on the objects of ${\cB}(\mC,\mA)$ the product

\begin{equation}
\label{15}
\phi \psi := \psi_{\rho,\rho} \circ \phi
\end{equation}

\noindent with identity the amplimorphism $\iota_\mA : \mA \ra (\iota,\iota)^\mA$ defined by $\iota_\mA (a) := a \otimes 1_\iota$. We now state some relations in order to define a tensor product on the arrows of $\cB(\mC,\mA)$. First of all note that the previous lemma implies that the linear map

$$
\psi_{\rho,\rho'}(x) := \psi(x_1)_1 \otimes \psi(x_1)_2 \times x_2 \in (\sigma \rho,\sigma \rho')^{\mA}
$$

\noindent is well defined as a Banach space map for $x \in (\rho,\rho')^\mA$ (note that the notation is consistent with (\ref{def_amp})). By verifying on the algebraic tensor products and as $\psi_{\rho,\rho'}$ is bounded we obtain, for $x' \in (\rho',\rho'')^\mA$

\begin{equation}
\label{eq_comp1}
\psi_{\rho,\rho''}(x'x) = \psi_{\rho',\rho''}(x') \cdot \psi_{\rho,\rho'}(x) .
\end{equation}

\noindent We also consider for each $\rho$ the Banach space map

$$
i_{\rho}(y) := y_1 \otimes y_2 \times 1_{\rho} \in ( \sigma \rho,\sigma' \rho)^{\mA}
$$

\noindent with $y \in (\sigma,\sigma')^{\mA}$. Let now $\phi,\phi',\psi,\psi'$ be amplimorphisms with $x \in (\rho,\rho')^{\mA} \subset (\phi,\phi')$, $y \in (\sigma,\sigma')^{\mA} \subset (\psi,\psi')$. Then for $t \in (\rho,\rho')^{\mA}$ we find, with the same argument of (\ref{eq_comp1}),

\begin{equation}
\label{eq_comp2}
i_{\rho'}(y) \cdot \psi_{\rho,\rho'} (t) = \psi_{\rho,\rho'}' (t) \cdot i_{\rho}(y) .
\end{equation}

\noindent So we define

\begin{equation}
\label{16}
x \times y := \psi_{\rho,\rho'}'(x) \cdot i_\rho (y) = i_{\rho'}(y) \cdot \psi_{\rho,\rho'}(x).
\end{equation}

\begin{prop}
\label{prop_tens}
Let $\mA$ be a unital \sC algebra and $\mC$ a strict tensor \sC category such that for each object $\rho $ of $\mC$ the maximal norm on $(\rho,\rho)_a^\mA$ coincides with the spatial norm. Then, endowed with the tensor product defined in (\ref{15}) and (\ref{16}), ${\cM}(\mC,\mA)$ is a strict semitensor \sC category and ${\cB}(\mC,\mA) $ is a strict tensor \sC category. 
\end{prop}

\begin{proof}
We verify that $x \times y \in (\phi \psi , \phi' \psi')$; by (\ref{eq_comp1}),(\ref{eq_comp2}) we obtain, with $a \in \mA$,

$$
\begin{array}{ll} 
x \times y \cdot \phi \psi (a) 

& = \psi_{\rho,\rho'}'(x) \cdot i_\rho (y) \cdot \psi_{\rho,\rho} (\phi (a)) \\
& = \psi_{\rho,\rho'}'(x) \cdot \psi_{\rho,\rho}' (\phi (a)) \cdot i_\rho (y) \\
& = \psi_{\rho,\rho'}'(x \phi (a)) \cdot i_\rho (y) \\
& = \psi_{\rho,\rho'}'(\phi' (a) x) \cdot i_\rho (y) \\
& = \psi_{\rho',\rho'}'(\phi' (a)) \cdot \psi_{\rho,\rho'}' (x) \cdot i_\rho (y) \\
& = \psi' \phi' (a) \cdot  x \times y \quad . \\

\end{array}
$$

\noindent The same computation for $a = 1$ shows that $x \times y \in (\phi \psi , \phi' \psi')_\mM$ if $x \in (\phi,\phi')_\mM$, $y \in (\psi,\psi')_\mM$. The other properties of tensor product can be verified by routine computations.
\end{proof}

Note that the space of intertwiners of the identity object $\iota_{\mA}$ is isomorphic to $(\iota,\iota)^{\mA}$ in the category ${\cM}(\mC,\mA)$ and to $(\iota,\iota)^{\mZ}$ in the category ${\cB}(\mC,\mA)$.

\begin{ex}
\label{ex61}
The motivating examples come from the case $\mC = {\bf Hilb}$. Then ${\mC}^{\mZ}$, ${\cB}(\mC,\mA)$, ${\cM}(\mC,\mA)$, ${\mC}^\mA$ correspond respectively to the categories of vector bundles over the spectrum of $\mZ$, Hilbert $\mA$-$\mA$-bimodules with arrows commuting with left and right action (\cite {DR89} \S 1.7), Hilbert $\mA$-$\mA$-bimodules, Hilbert $\mA$-modules. Consider in fact an homomorphism $\phi : \mA \ra {\bM}_n \otimes \mA$. Then we get a projective Hilbert $\mA$-$\mA$-bimodule $\mX := \left\{ x \in \mA^{\oplus^n}: \phi(1) x = x \right\}$; $\mA$ acts on $\mX$ on the left by $a,x \mapsto \phi(a)x$ and on the right by $x,a \mapsto xa$. The Hilbert $\mZ$-bimodule $\mX^\phi := \left\{ x \in \mX : \phi(a) x = xa, a \in \mA \right\}$ is also defined. As expected we have a tensor product in the categories of vector bundles and Hilbert $\mA$-$\mA$-bimodules (choosing the morphisms commuting with the actions), while only a semitensor structure is allowed on Hilbert $\mA$-$\mA$-bimodules, as well known.
\end{ex}

\begin{ex} Let $\mB$ be a unital \sC algebra with trivial center. The objects of ${\cB}({\bf end}\mB,\mA) $ are given by endomorphisms

$$
\phi : \mA \ra (\rho,\rho)^\mA \subset \mA \otimes \mB 
$$

\noindent where $\rho $ is an endomorphism on $\mB$. The identity object in ${\cB}({\bf end}\mB,\mA)$ is the natural immersion $\iota : \mA \hra (\iota_\mB,\iota_\mB)^\mA \simeq \mA \otimes 1_\mB$, where $\iota_\mB$ is the identity automorphism on $\mB$; this implies

$$
(\iota,\iota)_{\mM} = (\iota_\mB,\iota_\mB)^\mA \simeq \mA, 
$$

\noindent while

$$
(\iota,\iota) \simeq \left\{ x\in ( \iota_\mB,\iota_\mB)^\mA : a \otimes 1_\mB \cdot x=a\otimes 1_\mB \cdot x,a\in \mA\right\} \simeq \mZ. 
$$

\noindent Now for each $\phi$ the following Banach space is naturally defined:

$$
\mX_\phi := (\iota,\phi)_{\mM} := \left\{ x \in ( \iota_\mB,\rho ) ^{\mA}:\phi ( 1) x=x\right\}.
$$

\noindent $\mX_\phi $ is a Hilbert $\mA$-module in $\mA \otimes \mB$ with the usual $\mA$-valued scalar product and right multiplication. There is a left action of $(\phi,\phi)_{\mM}$ and this implies that a left action $a,x \mapsto \phi (a) x$ of $\mA$ on $\mX_\phi$ is defined.
\end{ex}

\bigskip

We now state some elementary properties of the DR-algebras associated to objects belonging to the categories appearing in (\ref{14}). We denote by the notation $\mO_{\phi,{\mM}}$ the DR-algebra associated to an amplimorphism $\phi$ as an object of $\cM(\mC,\mA)$. The following proposition is a categorical version of \cite{DPZ97} Prop. 3.4.

\begin{prop}
The inclusion functor ${\cB}(\mC,\mA) \stackrel{i} \hra {\cM}(\mC,\mA)$ induces a monomorphism $i: \mO_\phi \hra \mO_{\phi,\mM}$ such that $i(\mO_\phi) \subseteq \mA' \cap \mO_{\phi,{\mM}}$. If for some $ r \in \bN$ the space of arrows $(\iota_\mA,\phi^r)$ contains an isometry then $i(\mO_\phi) = \mA' \cap {\mO}_{\phi,{\mM}}$. 
\end{prop}

\begin{proof}
The proof follows the same line as in \cite{DPZ97}. For reader's convenience we expose some elementary facts about the first assertion. The embedding of the \sC algebras follows by funtoriality. Now, $\mA$ is embedded in $\mO_{\phi ,{\mM}}$ as the space of arrows of the identity object and if $t \in (\phi^r,\phi^{r+k})$ the multiplication by $a \in \mA$ is given by 

$$
\begin{array}{ll} 
at & = a \times \phi(1) \times ... \times \phi(1) \cdot t \\
& = a \times \phi^{r+k}(1) \cdot t \\  
& = \phi^{r+k}(a) t = t \phi^r(a) \\  
& = t \cdot a \times \phi^r(1) = ta
\end{array}
$$

\noindent and this proves that arrows between tensor powers of $\phi$ commute with elements of $\mA$. 
\end{proof}

The following results characterize DR-algebras in $\cM (\mC,\mA)$, $\cB (\mC,\mA)$ coming from 'commutative objects' belonging to ${\mC}^{\mZ}$.

\begin{lem}
\label{lem64}
Let $\mA$ be a unital \sC algebra, $\mC$ a strict tensor \sC category with identity object $\iota$ such that $(\iota,\iota) = \bC$, $e \in (\rho,\rho)^{\mZ}$ an object of ${\mC}^{\mZ}$ such that $(\iota,e)$ is finitely generated as a right Hilbert $\mZ$-module. If $\phi_e$ is the amplimorphism defined by $e$, then for $r+k \geq 0$:

\begin{itemize}

\item [{\em (1)} \quad]  $(e^r,e^{r+k}) = ( \phi_e^r , \phi_e^{r+k})$;

\item [{\em (2)} \quad]  $(\phi_e^r,\phi_e^{r+k})_{\mM}$ and $(e^r,e^{r+k}) \otimes_{\mZ} \mA$ are isomorphics as Banach $\mA$-$\mA$-bimodules.

\end{itemize}
\end{lem}

\begin{proof}
Let $t \in (e^r,e^{r+k})$. By definition of $\phi_e$ it follows that $t$ (we write $t$
for $i(t)$) belongs to $( \phi_e^r,\phi_e^{r+k}) $ if and only if commutes with $a\in \mA$; but

$$
\begin{array}{ll}
at & = a\times e^{r+k}\cdot t \\  
& = a( e^{r+k}) _1t_1\otimes (e^{r+k})_2t_2 \\
& = at_1\otimes t_2 = t_1a\otimes t_2\\
& = ta.
\end{array}
$$

\noindent Let now $t \in (\phi_e^r,\phi_e^{r+k})$; as $(\iota,e)$ is finitely
generated there exists a finite set $\left\{\psi_i\right\} \subset (\iota,e)$ such
that $e = \sum_i \psi_i \psi_i^*$ and this implies that each $(\iota ,e^r)$ is finitely generated by elements of the type $\psi_I := \psi_{i_1}...\psi_{i_r}$. Now we have seen that $( \iota ,e^r) \subset ( \iota ,\phi _e^r) $ so that $t_{IJ}:=\psi _I^{*}t\psi _J\in ( \iota ,\iota ) \simeq {\mZ}$ and $t=e^{r+k}t=te^r=\sum_{IJ}t_{IJ}\psi _I\psi _J^{*}$ belongs to $( e^r,e^{r+k}) $. For the second point, we define the map from $( e^r,e^{r+k}) \otimes_{\mZ} \mA$ into $( \phi _e^r,\phi _e^{r+k})_{\mM}$ by $t\otimes a:=at$. As $(at)^* a' t' = (a^* a') (t^* t')$ the map is injective, and by the argument of (1) if $x \in (\phi_e^r,\phi_e^{r+k})_{\mM}$ we get a decomposition $ x=\sum_{IJ}x_{IJ}\psi _I\psi _J^{*}$, $x_{IJ} \in \mA$, so we have the surjectivity.
\end{proof}

\begin{cor}
\label{cor65}
With the hypothesis of the previous lemma, there is an isomorphism of \sC dynamical systems $(\oex,\sigma_e) \simeq ( \mO_{\phi_e},\sigma_{\phi_e})$.
\end{cor}

We apply the previous results to Hilbert \sC bimodules. Let $\mX$ be the Hilbert \sC bimodule defined by the amplimorphism $\phi : \mA \ra \mA \otimes \bM_n$ as in Example \ref{ex61}. By using the notations of \cite {DPZ97} and referring to example \ref{ex61}, we get an isomorphism $\mO_{\phi,{\mM}} \simeq \mO_{\mX_{\mA}}$, where $\mO_{\mX_{\mA}}$ is the Pimsner \sC algebra associated to $\mX$, while $\mO_\phi \simeq \mO_{_{\mA}\mX_{\mA}}$, where $\mO_{_{\mA}\mX_{\mA}}$ is the subalgebra of $\mO_{\mX_{\mA}}$ generated by the arrows for which left and right actions of $\mA$ coincide. We denote by $(\mX^r , \mX^{r+k})$ the spaces of right $\mA$-module maps from $\mX^r$ into $\mX^{r+k}$, that are the arrows generating $\mO_{\mX_{\mA}}$. As we stated in the general case, there is an embedding $\mO_{_{\mA} \mX_{\mA}} \hra \mO_{\mX_{\mA}}$ such that $\mO_{_{\mA}\mX_{\mA}} \subseteq \mA^{\prime} \cap \mO_{\mX_{\mA}}$. We say that $\mA$ is {\em normal} in $\mO_{\mX_{\mA}}$ if $\mA =( \mA^{\prime} \cap \mO_{\mX_\mA})^{\prime} \cap \mO_{\mX_{\mA}}$ . Let now $\sigma_{\mX}$ be the canonical endomorphism defined on $\mO_{_{\mA}\mX_\mA}$. We say that $\sigma_\mX$ is {\em inner} if $\mX^\phi := \mX \cap \mO _{_\mA\mX_\mA}$ is finitely generated as Hilbert $\mZ$-module and with support $1$ in $\mO_{\mX_\mA}$.

\begin{prop}
\label{prop_mod}
Let $\mX$ be a Hilbert $\mA$-$\mA$-bimodule defined by an amplimorphism $\phi :\mA \ra {\bM}_n \otimes \mA$, $\mA$ normal in {\it $\mO_{\mX_\mA}$}. Then the following are equivalent:

\begin{itemize}
\item [{\em (1)} \quad]  There exists a vector bundle $\mE$ over the spectrum of $\mZ$ with module of section isomorphic to $\mX^\phi$ and $\mX \simeq \mX^\phi \otimes_{\mZ} \mA$;

\item [{\em (2)} \quad]  $\sigma_\mX$ is inner;

\item [{\em (3)} \quad]  $(\coe,\sigma_{\mE}) \simeq ( \mO_{_\mA \mX_\mA},\sigma_\mX )$ and $\mO_{{\mX}_\mA} \simeq \coe \otimes_{\mZ} \mA$

\end{itemize}

\end{prop}

\begin{proof}
As $\mX^\phi$ is finitely generated it is isomorphic to a projective Hilbert $\mZ$-module, hence isomorphic to the module of sections of a vector bundle $\mE$. Now $\sigma_\mX$ is inner if and only if the injective $\mA$-module map $\psi \otimes a := \psi a$ from $\mX^\phi \otimes_{\mZ} \mA$ to $\mX$ is surjective, so that the definition of the left action $\phi $ is forced by $\phi (a') (\psi a) := \psi a' a$. The equivalence with the third statement is consequence of Lemma \ref{lem64} and Corollary \ref{cor65} applied to the projection $e \in {\bM}_n \otimes \mZ$ defined by $\mX^\phi$: we obtain in this way the isomorphism $(\coe,\sigma_{\mE}) \simeq ( \mO_{_\mA \mX_\mA},\sigma_\mX )$, and a family of compatible Banach $\mA$-$\mA$-bimodule isomorphisms $(\mX^r,\mX^{r+k}) \simeq (\mE^r,\mE^{r+k}) \otimes_{\mZ} \mA$. This implies the isomorphism between the algebraic tensor product $^0\mO_\mE \otimes_\mZ \mA$ and $^0\mO_{{\mX}_\mA}$. Now observe that by nuclearity of $\coe$ there is a unique \sC norm defining the tensor product $\coe \otimes_{\mZ} \mA$ (in the sense of \cite{Kas88}, see also \cite{Bla00}), that coincides with the norm on $\mO_{{\mX}_\mA}$ for which the circle action is isometric.
\end{proof}

Let $\mE$ be a vector bundle over a compact Haussdorff space $\Omega$. We endow $\coe^1$ with a structure of a Hilbert $\coe^0$-$\coe^0$-bimodule by defining the right action $\psi \cdot x := \psi x$ and left action $x \cdot \psi := \sigma_{\mE} (x) \psi$, where $x \in \coe^0$, $\psi \in \coe^1$ and $\sigma_{\mE}$ is the canonical endomorphism on $\coe$. Note that $\coe^1$ is generated by $(\iota,\mE)$, on which left and right actions coincide. The previous proposition then implies (see \cite {Pim97} for an anologue result):

\begin{cor} Let $\mE$ be a vector bundle over a compact Haussdorff space $\Omega$. Then 

\begin{itemize}
\item [{\em (1)} \quad] $\coe^1 \simeq ( \iota ,{\mE}) \otimes _{C(\Omega)} \coe^0$

\item [{\em (2)} \quad] {\it $\mO_{\mO_{\mE}^1}$}$ \simeq \coe \otimes_{C(\Omega)} \coe^0.$

\end{itemize}

\end{cor}

\section{Final Remarks.}

This paper does not claim to be exhaustive. There are several interesting open questions that are the subject of the previously quoted work in progress between S. Doplicher, P. Goldstein and the author, first of all the geometrical characterization of the class of vector bundles having isomorphic DR-algebras. Another open question is the computation of the $K$-theory of such \sC algebras; note that for DR-algebras \cite{Exe94,AEE95} a Pimsner-Voiculescu exact sequence holds, so the first step should be the computation of the $K$-theory of the zero grade algebra, that is an inductive limit of continuous fields of matrix algebras.

\bigskip

The notion of extension can be used to introduce $K$-type functors. Given in fact a \sC category $\mC$ with direct sums and a \sC algebra $\mA$ we can consider the set of unitary equivalence classes of objects in $\mC^\mA$; we endowe it with the operation of direct sum (that is defined on $\mC^\mA$ by Prop. \ref{prop_dir_sum}) and introduce the associated Grothendieck group $K_\mC (A)$, that coincides with the usual $K$-theory of $\mA$ in the case $\mC = {\bf Hilb}$. Properties as funtoriality, homotopic invariance and additivity can be easily verified. We defer to a future work the analysis of the other homology properties (continuity, half exactness), as well as the study of interesting cases (as the closure for direct sums of $\mC = {\bf end}\mB$ in Example \ref{enda}) .

\bigskip

A further interesting point is the extension of Doplicher-Roberts duality for compact groups to categories with not trivial space of intertwiners of the identity object. An important step in the theory is the crossed product of a \sC algebra by an endomorphism satisfying certain properties (namely permutation simmetry and special conjugate property, see \cite{DR89A}); in that case the \sC algebra has trivial center and the endomorphism is implemented by a Hilbert space in the cross product, while in the cases discussed here (for example the \sC algebra of a vector bundle) the canonical endomorphism is implemented by a finitely generated module over the center. The construction of such a generalized crossed product of a \sC algebra with not trivial center by an endomorphism should give the first step towards an extension of Doplicher-Roberts duality in the sense stated above. Some work in this direction has been done in \cite{BL97}, in the setting of Hilbert \sC systems. In this context free right modules over the center of the given \sC algebra appear, but the example of the \sC algebra of a vector bundle shows that more general objects can be considered, as the (not free) modules of sections of vector bundles.

\bigskip

{\bf Acknowledgments.} I'm very grateful to Professor S. Doplicher for having proposed the initial idea and for helpful suggestions and comments, and indebted to F. Lled\'o, G. Ruzzi, M. Gabriel for valuable discussions, interesting remarks and precious help.

\end{document}